%% file: quot3_new.tex
\documentclass[a4j,10pt]{article}

\usepackage{amsmath,amscd,amssymb}

\input{ntn}
\input{classdef}

\begin{document}

\title{
The structure of the cohomology ring
of the filt schemes
}
\author{Takuro Mochizuki}
\date{}

\maketitle

\abstract{
Let $C$ be a smooth projective curve
over the complex number field $\cnum$.
We investigate the structure of the cohomology ring
of the quot schemes $\Quot(r,n)$, i.e.,
the moduli scheme of the quotient sheaves of
$\nbigo_C^{\oplus\, r}$ with length $n$.
We obtain a filtration on $H^{\ast}(\Quot(r,n))$,
whose associated graded ring has a quite simple structure.
As a corollary, we obtain a small generator of the ring.
We also obtain a precise combinatorial description
of $H^{\ast}(\Quot(r,n))$ itself.

For that purpose,
we consider the complete filt schemes and
we use a `splitting principle'.
A complete filt scheme has an easy geometric description:
a sequence of bundles of projective spaces.
Thus the cohomology ring of the complete filt schemes are
easily determined.
In the cohomology ring of complete filt schemes,
we can do some calculations
for the cohomology ring of quot schemes.


}

\input{12.5}


\input{12.0}

\input{12.1}

\input{12.2}



\input{12.3}

\input{12.4}

\input{12.42}

\input{12.41}


\input{12.6}


\input{quotref}
\mbox{}{ }{ }{ }

\noindent
{\it Address\\
Department of Mathematics,
Osaka City University,
Sugimoto, Sumiyoshi-ku,
Osaka 558-8585, Japan.
takuro@sci.osaka-cu.ac.jp\\
School of Mathematics,
Institute for Advanced Study,
Einstein Drive, Princeton,
NJ, 08540, USA. \\
takuro@math.ias.edu
}

\end{document}

%% file: ntn.tex
\newcommand{\nbiga}{\mathcal{A}}
\newcommand{\nbigb}{\mathcal{B}}

\newcommand{\nbigd}{\mathcal{D}}

\newcommand{\nbigf}{\mathcal{F}}
\newcommand{\nbigg}{\mathcal{G}}

\newcommand{\nbigi}{\mathcal{I}}
\newcommand{\nbigj}{\mathcal{J}}
\newcommand{\nbigk}{\mathcal{K}}
\newcommand{\nbigl}{\mathcal{L}}

\newcommand{\nbign}{\mathcal{N}}
\newcommand{\nbigo}{\mathcal{O}}
\newcommand{\nbigp}{\mathcal{P}}

\newcommand{\nbigt}{\mathcal{T}}
\newcommand{\nbigu}{\mathcal{U}}

\newcommand{\proj}{\Bbb{P}}
\newcommand{\seisuu}{\Bbb{Z}}
\newcommand{\seisuuplus}{\Bbb{Z}_{\geq\,0}}
\newcommand{\rnum}{\Bbb{Q}}

\newcommand{\cnum}{{\Bbb C}}

\newcommand{\gbigs}{\frak S}

\newcommand{\vece}{{\boldsymbol e}}

\newcommand{\vecv}{{\boldsymbol v}}
\newcommand{\vecu}{{\boldsymbol u}}
\newcommand{\vecw}{{\boldsymbol w}}

\newcommand{\vecl}{{\boldsymbol l}}

\newcommand{\vecc}{{\boldsymbol c}}
\newcommand{\vech}{{\boldsymbol h}}

\newcommand{\vecL}{{\boldsymbol L}}

\newcommand{\rarr}{\rightarrow}
\newcommand{\lrarr}{\longrightarrow}

\newcommand{\pf}{{\bf Proof}\hspace{.1in}}
\newcommand{\qed}{\mbox{\rule{1.2mm}{3mm}}}

\def\Cok{\mathop{\rm Cok}\nolimits}
\def\cok{\mathop{\rm Cok}\nolimits}
\def\Image{\mathop{\rm Im}\nolimits}

\def\Gr{\mathop{\rm Gr}\nolimits}

\def\rank{\mathop{\rm rank}\nolimits}

\def\ker{\mathop{\rm Ker}\nolimits}

\def\Gr{\mathop{\rm Gr}\nolimits}
\def\Sym{\mathop{\rm Sym}\nolimits}
\def\co{\mathop{\rm co}\nolimits}
\def\Nor{\mathop{\rm Nor}\nolimits}

\newcommand{\Dec}{{\mathcal Dec}}

\def\Quot{\mathop{\rm Quot}\nolimits}
\def\Filt{\mathop{\rm Filt}\nolimits}
\def\res{\mathop{\rm res}\nolimits}
\def\sym{\mathop{\rm sym}\nolimits}

\newcommand{\projnbigku}{\proj(\nbigk_n^{u\,\lor}) }

\newcommand{\nbigfu}{\nbigf^u}
\newcommand{\nbiggu}{\nbigg^u}
\newcommand{\nbigku}{\nbigk^u}

\newcommand{\pt}{pt}
\newcommand{\leqzero}{\leq_0}
\newcommand{\geqzero}{\geq_0}
\newcommand{\leqiti}{\leq_1}

\newcommand{\shat}{\hat{S}}

%% file: classdef.tex
\newcommand{\bdmath}{\begin{displaymath}}
\newcommand{\edmath}{\end{displaymath}}

\newcommand{\beqn}{\begin{equation}}
\newcommand{\eeqn}{\end{equation}}

\newcommand{\beqnarray}{\begin{eqnarray}}
\newcommand{\eeqnarray}{\end{eqnarray}}

\newcommand{\bitemize}{\begin{itemize}}
\newcommand{\eitemize}{\end{itemize}}

\newcommand{\benumerate}{\begin{enumerate}}
\newcommand{\eenumerate}{\end{enumerate}}

\newcommand{\bdescriprion}{\begin{description}}
\newcommand{\edescriprion}{\end{description}}

\newtheorem{thm}{Theorem}[section]
\newtheorem{cor}{Corollary}[section]

\newtheorem{rem}{Remark}[section]
\newtheorem{lem}{Lemma}[section]
\newtheorem{prop}{Proposition}[section]
\newtheorem{df}{Definition}[section]

\def\longuparrow{\vcenter{%
   \lineskip0pt\lineskiplimit0pt\baselineskip0pt\ialign{%
   \hfil{##}\hfil\crcr\hbox to 0pt{\hss$\uparrow$\hss}\cr%
   \hbox to 0pt{\hss\vrule width.4pt depth 0pt height 1em\hss}\cr}}}
\def\longdownarrow{\vcenter{%
   \lineskip0pt\lineskiplimit0pt\baselineskip0pt\ialign{\hfil{##}\hfil%
   \crcr\hbox to 0pt{\hss\vrule width.4pt depth 0pt height 1em\hss}\cr%
   \hbox to 0pt{\hss$\downarrow$\hss}\cr}}}

%% file: 12.5.tex
\section{Introduction}

Let $C$ be a smooth projective curve over the complex number field
$\cnum$.
Let $E$ be a locally free coherent sheaf over $C$.
Let $n$ be a non-negative integer.
Let $\Quot(E,n)$ denote the moduli scheme of quotient sheaves
of $E$ with length $n$.
It is called the quot scheme.
In particular, we use the notation $\Quot(r,n)$
if $E=\nbigo^{\oplus\,r}$.
In the papers \cite{m1} and \cite{m2},
we studied the structure of the cohomology ring
of quot schemes $\Quot(r,l)\,\,(l=2,3)$.
We also studied the infinite quot schemes
$\bigcup_{r=1}^{\infty} \Quot(r,l)=\Quot(\infty,l)$.
In particular, we determined the generators and 
the relations among them
for the cohomology rings
$H^{\ast}(\Quot(\infty,l)),\,\,H^{\ast}(\Quot(r,l))\,\,(l=2,3)$.
However, when  $l\geq 4$,
the rings seem somewhat complicated for a direct calculation.
Changing the route for an attack,
we introduce the `splitting principle'
for the quot schemes,
which is an analogue of
the well known splitting principle in the theory of Grassmannian manifolds.

We introduce the filt schemes:
Let $E$ be a vector bundle over $C$,
and $\vecl=(l_1,\ldots,l_h)$ be a tuple of non-negative integers.
The filt scheme $\Filt(E,\vecl)$ is the moduli of
the sequence of quotients:
\[
 E\lrarr \nbigf_h\lrarr \nbigf_{h-1}\lrarr \cdots\lrarr
  \nbigf_2\lrarr\nbigf_1.
\]
It satisfies the following condition:
\begin{itemize}
\item
 We put $\nbigg_i:=\ker(\nbigf_i\lrarr \nbigf_{i-1})$.
Then the length of $\nbigg_i$ is $l_i$.
\end{itemize}
If $l_1=n$ and $l_i=0$ for any $i>1$,
the filt scheme $\Filt(E,(n,0,\ldots,0))$ is naturally identified with
the quot scheme $\Quot(E,n)$.
On the other hand,
if $l_i=1$ for any $i$,
the variety $\Filt(E,(\overbrace{1,\ldots,1}^n))$ is called
the complete filt scheme.
We denote it by $\Filt(E,n)$ for simplicity of the notation.
In particular, we denote $\Filt(\nbigo^{\oplus\,r},n)$
by $\Filt(r,n)$.
We put $\Filt(\infty,n):=\bigcup_{r=1}^{\infty}\Filt(r,n)$.

The structure of $\Filt(E,n)$ is quite simple.
As is described in the subsection \ref{subsection;1.9.80},
it is a projective space bundle
over $\Filt(E,n-1)\times C$.
Let $\nbigl_n$ denote the corresponding tautological line bundle.
We denote the first Chern class of $\nbigl_n$
by $\omega_n$.
We have the natural morphism
$\eta_{n,m}:\Filt(E,n)\lrarr \Filt(E,m)$ for any $m\leq n$,
by forgetting the parts
$\nbigf_{j}$ $(m+1\leq j\leq n)$.
We denote the line bundle $\eta_{n,m}^{\ast}\nbigl_m$
and $\eta_{n,m}^{\ast}\omega_m$
by $\nbigl_m$ and $\omega_m$ respectively.
We have the naturally defined morphism
$\Filt(E,n)\lrarr C^n$ by taking the supports
of $\nbigg_i$.
Then the structure of the cohomology ring $H^{\ast}(\Filt(E,n))$
over $H^{\ast}(C^n)$
is easily and completely described
in termes of $\omega_i$ $(i=1,\ldots,n)$.
(See the subsubsection \ref{subsubsection;1.9.81}.)
When we take a limit,
$H^{\ast}(\Filt(\infty,n))$ is
the polynomial ring over $H^{\ast}(C^n)$
generated by $\omega_i$ $(i=1,\ldots,n)$.

We have the natural morphism
$\Psi:\Filt(E,n)\lrarr \Quot(E,n)$,
forgetting $\nbigf_i$ $(i=1,\ldots,n-1)$.
It induces the injection of the cohomology rings.
Thus $H^{\ast}(\Quot(E,n))$
is a subring of the ring
$H^{\ast}(\Filt(E,n))$.
We would like to do some calculation
for $H^{\ast}(\Quot(E,n))$
in $H^{\ast}(\Filt(E,n))$.

Let consider the case $E=\bigoplus_{\alpha=0}^{r-1}L_{\alpha}$.
Let $G_m^r$ denote the $r$-dimensional algebraic torus.
We have the decomposition of the fixed point sets:
\[
 \Quot(E,n)^{G_m^r}
=\coprod_{\vecv\in \nbigb(n,r)}
 F^Q(\vecv),
\quad
 \Filt(E,n)^{G_m^r}
=\coprod_{\vecv\in[0,r-1]^n}
 F(\vecv).
\]
(See the subsubsection \ref{subsubsection;1.9.85}.)
Due to the theories of Bianicki-Birula (\cite{bb}),
Carrel and Sommes (\cite{cs}),
or in other words,
due to the Morse theoretic consideration,
we obtain the morphisms:
\[
\begin{array}{ll}
 \xi^Q(\vecv;\cdot):
 H^{\ast}(F^Q(\vecv))[-2\co(\vecv)]
 \lrarr H^{\ast}(\Quot(E,n)),
 &(\vecv\in\nbigb(n,r))
 \\
 \mbox{{}}\\
 \xi(\vecv;\cdot):
 H^{\ast}(F(\vecv))[-2\co(\vecv)]
 \lrarr H^{\ast}(\Filt(E,n)),
 &(\vecv\in [0,r-1]^n).
\end{array} 
\]
And we have the decompositions:
\[
 H^{\ast}(\Quot(E,n))
=\bigoplus_{\vecv\in\nbigb(n,r)}
 \Image(\xi^Q(\vecv;\cdot)),
\quad
  H^{\ast}(\Filt(E,n))
=\bigoplus_{\vecv\in[0,r-1]^n}
 \Image(\xi(\vecv;\cdot))
\]
We have the following relation (the subsection \ref{subsection;1.13.1}):
\begin{equation} \label{eq;1.9.90}
 \Psi^{\ast}(\xi^Q(\vecv,a))=
\frac{1}{|St(\vecv)|}
\sum_{\sigma\in\gbigs_n}
 \xi(\sigma \vecv,\sigma(a)).
\end{equation}
Thus the calculation in  $H^{\ast}(\Quot(r,n))$
is reduced to the calculation in $H^{\ast}(\Filt(r,n))$.
The main calculation will be done in 
Theorem \ref{thm;1.9.40}.
(See the formula (\ref{eq;1.9.71}).
 See also (\ref{eq;1.9.70}).)
Although the equality (\ref{eq;1.9.71}) provides us
the method of the calculation,
the explicit calculation seems rather complicated.
Thus we introduce the filtration defined as follows:
\[
 F_h H^{\ast}(\Quot(E,n))
=\bigoplus_{
 \substack{\vecv\in\nbigb(n,r)\\ \co(\vecv)\leq h}}
 \Image(\xi^Q(\vecv;\cdot)),
\quad
 F_h H^{\ast}(\Filt(E,n))
=\bigoplus_{\substack{\vecv\in[0,r-1]^n\\ \co(\vecv)\leq h}}
 \Image(\xi(\vecv;\cdot)).
\]
Due to the equality (\ref{eq;1.9.90})
and Theorem \ref{thm;1.9.40},
it is easy to derive that the filtrations
are compatible with the products.
Moreover the structures of the associated graded rings
are quite simple (the subsubsectoin \ref{subsubsection;1.13.5}).

As a corollary,
we obtain some sets of the generators of 
$H^{\ast}(\Quot(\infty,n))$.
(Proposition \ref{prop;1.13.10}, Proposition \ref{prop;1.13.11}
and Theorem \ref{thm;1.13.12}.)

To determine the precise structure of $H^{\ast}(\Quot(r,n))$,
we consider the $\gbigs_n$-action $\rho$
on $H^{\ast}(\Filt(\infty,n))=H^{\ast}(C^n)[\omega_1,\ldots,\omega_n]$,
defined as follows (See subsection \ref{subsection;1.14.20}):
\begin{quote}
For $\sigma\in\gbigs_n$,
we put
$\rho(\sigma)(\prod_i p_i^{\ast}(a_i)\cdot\omega_i^{l_i})
  :=\prod_i p_{\sigma(i)}^{\ast}(a_i)\cdot\omega_{\sigma(i)}^{l_i}
 $.
\end{quote}
The action $\rho$ preserves the products.
Let $H^{\ast}(C^n)[\omega_1,\ldots,\omega_n]^{\rho,\gbigs_n}$
denote the $\gbigs_n$-invariant part with respect to the action
$\rho$.
Then we have
$\Psi^{\ast}H^{\ast}(\Quot(\infty,n))=
 H^{\ast}(C^n)[\omega_1,\ldots,\omega_n]^{\rho,\gbigs_n}$.
(Theorem \ref{thm;1.14.21})

\begin{rem}
We also use the $\gbigs_n$-action $\kappa$
on $H^{\ast}(\Filt(r,n))$ defined by
$\kappa(\sigma)\xi(\vecv,a)=\xi(\sigma\cdot\vecv,\sigma(a))$.
(See the subsection {\rm \ref{subsection;1.14.1}}.
We omit to denote $\kappa$ there.)
We also have
$H^{\ast}(\Filt(r,n))^{\kappa,\gbigs_n}=
 \Psi^{\ast}H^{\ast}(\Quot(r,n))$.
However the action $\kappa$ does not preserve the product.
\hfill\qed
\end{rem}

Let $\nbigj(r,n)$ denote the ideal of $H^{\ast}(\Quot(r,n))$
generated by 
$\{\xi^Q(\vecv)\,|\,\vecv\in\nbigb(n,r+1)-\nbigb(n,r)\}$.
Then we have
$H^{\ast}(\Quot(r,n))=H^{\ast}(\Quot(\infty,n))/\nbigj(r,n)$.
Thus, once we obtain
the description of $\Psi^{\ast}\xi^Q(\vecv;a)$
in terms of $\omega_i$,
we can say that 
the structure $H^{\ast}(\Quot(r,n))$ is described.

The formulas (\ref{eq;1.9.71}) and (\ref{eq;1.9.70})
assures us that we have some algorithms to calculate
$\Psi^{\ast}\xi^Q(\vecv;a)$.
Each algorithm will provide us a combinatorial decription.
In the subsection \ref{subsection;1.15.1},
we give one of the description.
Although the result does not seem convenient for the practical calculation,
we would like to remark the following:
Due to Theorem \ref{thm;1.14.21},
the formula gives a symmetric function,
that is, they are elements of
$H^{\ast}(C^n)[\omega_1,\ldots,\omega_n]^{\rho,\gbigs_n}$,
which is not obvious at a sight.

Finally we see the structure of the cohomology ring
of the infinite quot scheme of infinite length $\Quot(\infty,\infty)$.

\vspace{.1in}

The filt schemes $\Filt(E,\vecl)$ are the quite natural objects:
They are defined as the moduli scheme of some simple functors.
They are smooth.
They also have big symmetry.
The author believes that
it is always significant for mathematics
to investigate such objects in detail,
and he expects that the they will provide us
very interesting examples.

This paper is a revision.
The most main results are same as those of the original version.
The original version was written in May 2001,
after the original versions of the papers \cite{m1} and \cite{m2}.
The subsections \ref{subsection;1.14.20} and \ref{subsection;1.15.1}
are added in January 2003.

\vspace{.1in}

\noindent
{\bf Acknowledgement}\\
The author is grateful to the colleagues in Osaka City University.
In particular, he thanks Mikiya Masuda for his encouragement and support.

The author thanks the financial supports by Japan Society for
the Promotion of Science and the Sumitomo Foundation.
The paper was revised during his stay at the Institute
for Advanced Study. The author is sincerely grateful
to their excellent hospitality.
He also acknowledges National Scientific Foundation
for a grant DMS 9729992,
although any opinions, findings and conclusions or recommendations
expressed in this material are those of the author.

%% file: 12.0.tex
\section{Preliminary}

\subsection{Notation}
\label{subsection;1.13.30}
We use the following notation:
\begin{center}
\begin{tabular}{lll}
 $\seisuu$: The ring of integers, &
 $\seisuuplus:=\{a\in\seisuu\,|\,a\geq 0\}$,&
 $[0,n]:=\{a\in\seisuu\,|\,0\leq a\leq n\}$,\\
\mbox{{}}\\
 $[1,n]:=\{a\in\seisuu\,|\,1\leq a\leq n\}$,&
 $\cnum$: The complex number field, &
 $\rnum$: The rational number field.\\
\mbox{{}}\\
 $\gbigs_n$: The $n$-th symmetric group.
\end{tabular}
\end{center}

\begin{itemize}
\item
The order of a finite set $S$ is denoted by $|S|$.
\item
In general, $q_j$ will denote the projection
onto the $j$-th component.
For example, $q_j(\vecl)$ denotes the $j$-th component $l_j\in\seisuuplus$
of $\vecl=(l_1,l_2,\ldots,l_n)\in\seisuuplus^n$.
\item
For $1\leq j,k\leq n$, the transposition of $j$ and $k$
is denoted by $\tau_{j,k}$.
Namely $\tau_{j,k}$ denotes the element of $\gbigs_n$,
satisfying $\tau_{j,k}(j)=k$, $\tau_{j,k}(k)=j$
and $\tau_{j,k}(l)=l$ for any $l\neq j,k$.
\item
Let $n$ be a non-negative integer,
and $X$ be a variety.
Then $X^{(n)}$ denote the $n$-th symmetric product
of $X$.
For an element $\vecl=(l_1,\ldots,l_h)\in\seisuuplus^h$,
we put $X^{(\vecl)}:=\prod_{i=1}^h X^{(l_i)}$.

\item
The element
$(\overbrace{0,\ldots,0}^{i-1},1,0,\ldots,0)\in
 \seisuu_{\geq\,0}^n$
is denoted by $\vece_i$.
We put $\vecc_l:=\sum_{i=1}^l\vece_i$.

\item
For any element $\vecl\in\seisuuplus^n$,
we put $|\vecl|:=\sum_{i=1}^n q_i(\vecl)$.
Here $q_i(\vecl)$ denotes the $i$-th component of $\vecl$.
\item
In general, the identity morphism of $X$ is denoted by $id_X$.
We often omit to denote $X$, if there are no confusion.

\item
Let $C$ be a smooth projective curve.
We denote the diagonal of $C\times C$ by $\Delta$.
Let $\Delta_{i\,j}$ $(i\neq j)$ denote
the divisor
$\{(x_1,\ldots,x_n)\in C^n\,|\,x_i=x_j\}$ of $C^n$.
Let $I$ be a subset of $\{1,\ldots,n\}$,
then $\Delta_I$ denotes the closed subset
$\{(x_1,\ldots,x_n)\,|\,x_i=x_j\,\mbox{ for any $i,j\in I$}\}$.

We will often use the notation $\Delta_I$
to denote the corresponding cohomology class.

\item
In general,
we use the notation $\pt$ to denote the cohomology class
of a point in a topological space $X$.
Let $C$ be a smooth projective curve.
For a subset $I\subset \{1,\ldots,n\}$,
we have the projection $\prod_{i\in I}q_i:C^n\lrarr C^{|I|}$.
We have the cohomology class $\pt$ in $H^{\ast}(C^{|I|})$
in the above sense.
The pull back $(\prod_{i\in I}q_i)^{\ast}\pt$ is denoted by $\pt_I$.

\item
 In this paper, the coefficient of the cohomology ring
is the rational number field $\rnum$,
if we do not mention.
\item
We will consider the polynomial ring $H^{\ast}(X)[t_0,\ldots,t_{r-1}]$
over $H^{\ast}(X)$.
Here $X$ denotes a topological space, and $t_i$ $(i=0,\ldots,r-1)$
denote formal variables.
For an element $J=(j_0,\ldots,j_{r-1})\in \seisuuplus^r$,
we put $t^J:=\prod_{i=0}^{r-1}t_i^{j_i}$.
Then an element $P\in H^{\ast}(X)[t_0,\ldots,t_{r-1}]$
has the description $P=\sum_J a_J\cdot  t^J$
$\bigl(a_J\in H^{\ast}(X)\bigr)$.
The total degree of $P$ with respect to the variables
$t_0,\ldots,t_{r-1}$
is denoted by $\deg_t(P)$.
It is called total degree for simplicity.
We also put as follows:
\[
 [P]:=
 \sum_{|J|=\deg_t(P)}a_J\cdot t^J.
\]
It is called the top term of $P$.
\item
In principle, the polynomial ring $H^{\ast}(X)[t_0,\ldots,t_{r-1}]$
is obtained as the equivariant cohomology ring
$H^{\ast}_{G_m^r}(X)$ for the trivial action of $G_m^r$.
For an element $a\in H^{j}_{G_m^r}(X)$,
the number $j$ is called the cohomological degree of $a$.
\end{itemize}

\subsection{Some sets}

\subsubsection{
${\boldsymbol\nbigb\boldsymbol(\boldsymbol l\boldsymbol)}$
and
${\boldsymbol\nbigb\boldsymbol(\boldsymbol l,\boldsymbol r\boldsymbol)}$}
Let $l$ be a non-negative integer.
We put as follows:
\[
 \nbigb(l):=
 \bigl\{
 \vecv=(v_1,\ldots,v_l)\in\seisuuplus^l\,
 \big|\,
 v_1\geq v_2\geq \cdots\geq v_l\geq 0
 \bigr\}.
\]
Let $r$ be a positive integer.
Then we put as follows:
\[
 \nbigb(l,r):=
 \bigl\{
 \vecv=(v_1,\ldots,v_l)\in\nbigb(l)\,
 \big|\,
 v_1\leq r-1
 \bigr\}.
\]
We have the natural inclusion $\nbigb(l)\lrarr \seisuuplus^{l}$.
Then the set $\nbigb(l,r)$ is contained in $[0,r-1]^l$.

Let $\vecu=(u_1,\ldots,u_l)$ be any element of $\seisuuplus^l$.
We have the unique element $\vecv$ of $\nbigb(l)$,
which is a reordering of $\vecu$.
Thus we have the map $\seisuuplus^l\lrarr\nbigb(l)$,
which we denote by $\Nor$,
i.e. $\vecv:=\Nor(\vecu)$.
Clearly we have $\Nor([0,r-1]^l)\subset \nbigb(r,l)$.
The map is called the normalization.

For an element $\vecv\in \nbigb(l)$,
we put $\co(\vecv):=\sum_{i=1}^l v_i\geq 0$.
Thus we obtain the map
$\co:\nbigb(l)\lrarr \seisuuplus$.

For a tuple $\vecl=(l_1,\ldots,l_h)\in\seisuuplus^h$,
we put as follows:
\[
 \nbigb(\vecl):=\prod_{i=1}^h\nbigb(l_i),
\quad
 \nbigb(\vecl,r):=\prod_{i=1}^h\nbigb(l_i,r).
\]
For an element
$\vecv_{\ast}=(\vecv_1,\ldots,\vecv_h)
 \in \nbigb(\vecl)$,
we put
$ \co(\vecv_{\ast}):=\sum_{i=1}^h \co(\vecv_i)
\in\seisuuplus$.
For an element $\vecu\in\prod_{i=1}^h\seisuuplus^{l_i}$,
we put as follows:
\[
 \Nor(\vecu):=
 \bigl(
 \Nor(\vecu_1),\ldots,\Nor(\vecu_h)
 \bigr)
 \in \nbigb(\vecl).
\]
Thus we obtain the maps
$\co:\nbigb(\vecl)\lrarr\seisuuplus$
and
$\Nor:\prod_{i=1}^h\seisuuplus^{l_i}\lrarr \nbigb(\vecl)$.

\vspace{.1in}
\noindent
{\bf Example}\\
We have the natural isomorphism
$\nbigb(\vecc_n)\simeq \seisuuplus^n$.
We have the natural isomorphism
$\nbigb(n\cdot\vece_1)\simeq\nbigb(n)$.

\subsubsection{$\boldsymbol{\Dec(\vecl,r)}$}

Let $\vecl$ be an element of $\seisuuplus^h$.
A tuple
$\vecl_{\ast}=
 (\vecl_0,\ldots,\vecl_{r-1})
  \in 
 \bigl(
 \seisuuplus^{h}
 \bigr)^{r}$ is called an $r$-decomposition of $\vecl$,
if we have $\sum_{\alpha=0}^{r-1}\vecl_{\alpha}=\vecl$.
Let $\Dec(\vecl,r)$ denote the set of $r$-decompositions
of $\vecl$.
We have the natural injection
$\Dec(\vecl,r)\lrarr \Dec(\vecl,r+1)$
given by
$(\vecl_0,\ldots,\vecl_{r-1})\longmapsto 
 (\vecl_0,\ldots,\vecl_{r-1},0)$.
We put $\Dec(\vecl):=\bigcup_{r\geq 1}\Dec(\vecl,r)$.

For an element $\vecl_{\ast}=(\vecl_0,\vecl_1,\ldots,\vecl_{r-1})$,
we put as follows:
\[
 \co(\vecl_{\ast}):=
 \sum_{\alpha=0}^{r-1}
 \alpha\cdot |\vecl_{\alpha}|
=\sum_{\alpha=0}^{r-1}
 \alpha\cdot
 \sum_{i=1}^hq_i(\vecl_{\alpha}).
\]

\subsubsection{Isomorphism}

Let $\vecl$ be an element of $\seisuuplus^{h}$.
We have the natural bijection $\phi$ between
$\nbigb(\vecl,r)$ and $\Dec(\vecl,r)$.
Let $\vecv_{\ast}=(\vecv_1,\ldots,\vecv_h)$ be an element 
of $\nbigb(\vecl,r)$.
Then $\phi(\vecv_{\ast})=\vecl_{\ast}=(\vecl_0,\ldots,\vecl_{r-1})$
is determined as follows:
\[
 q_j(\vecl_{\alpha}):=
\bigl|
 \bigl\{
 i\,\big|\,
 q_i(\vecv_{j})=\alpha
 \bigr\}
\bigr|.
\]
More visually, the correspondence is as follows:
\[
 \begin{array}{c}
 (v_{1\,1},v_{1\,2},\ldots,v_{1\,l_1})=
 (
 \overbrace{r-1,\ldots,r-1}^{l_{r-1,1}},
 \overbrace{r-2,\ldots,r-2}^{l_{r-2,1}},
  \ldots\ldots,
 \overbrace{1,\ldots,1}^{l_{1,1}},
 \overbrace{0,\ldots,0}^{l_{0,1}},
)\\
 \mbox{{}}\\
 (v_{2\,1},v_{2\,2},\ldots,v_{2\,l_2})=
 (\overbrace{r-1,\ldots,r-1}^{l_{r-1,2}},
 \overbrace{r-2,\ldots,r-2}^{l_{r-2,2}},
 \ldots\ldots,
 \overbrace{1,\ldots,1}^{l_{1,2}},
 \overbrace{0,\ldots,0}^{l_{0,2}},
 )\\
\mbox{{}}\\
\vdots\\
\mbox{{}} \\
  (v_{n,1},v_{n,2},\ldots,v_{n,l_n})=
 (\overbrace{r-1,\ldots,r-1}^{l_{r-1,n}},
  \overbrace{r-2,\ldots,r-2}^{l_{r-2,n}},
\ldots\ldots,
  \overbrace{1,\ldots,1}^{l_{1,n}},
  \overbrace{0,\ldots,0}^{l_{0,n}},
 ).
 \end{array}
\]
By the construction,
we clearly have $\co(\phi(\vecv_{\ast}))=\co(\vecv_{\ast})$.

%% file: 12.1.tex
\subsection{Filt scheme}
\label{section;1.13.20}

\subsubsection{Definition}

Let $\vecl=(l_1,\ldots,l_h)$ be an element of $\seisuu^{h}_{\geq 0}$.
Let $C$ be a smooth projective curve over $\cnum$.
Let $E$ be a locally free coherent sheaf over $C$.
The filt scheme $\Filt(E,\vecl)$ is the moduli scheme of
the sequence of quotients
\begin{equation} \label{eq;1.8.1}
 E\lrarr \nbigf_h\lrarr \nbigf_{h-1}\lrarr \cdots\lrarr
  \nbigf_2\lrarr\nbigf_1,
\end{equation}
satisfying the following conditions:
\begin{itemize}
\item
The length of the sheaf
$\nbigg_i=Ker(\nbigf_i\lrarr \nbigf_{i-1})$
is $l_i$ for each $i$.
\end{itemize}

In particular, $\Filt(E,\vecc_n)$ is denoted by
$\Filt(E,n)$.
It is called the complete filt scheme.
On the other hand,
$\Filt(E,n\cdot\vece_1)$ is denoted by $\Quot(E,n)$.
It is the usual quot scheme.

\subsubsection{Deformation}

A sequence as in (\ref{eq;1.8.1})
is denoted by $E\rarr \nbigf_{\ast}$.
For such sequence,
we put $\nbigk_j:=\ker(E\lrarr \nbigf_j)$
and $\nbigk_0:=E$.
Then we have the natural morphisms:
\[
 Hom(\nbigk_i,\nbigk_{j})\lrarr
 Hom(\nbigk_i,\nbigk_{j-1}),
\quad
 Hom(\nbigk_i,\nbigk_j)\lrarr
 Hom(\nbigk_{i+1},\nbigk_j).
\]
In all, we obtain the following morphism $f$:
\[
\begin{CD}
 \bigoplus_{i=1}^h
 Hom(\nbigk_i,\nbigk_{i})
@>{f}>>
 \bigoplus_{i=1}^h
 Hom(\nbigk_{i},\nbigk_{i-1}).
\end{CD}
\]
Here the first term stands at degree $-1$.
It is easily checked that $f$ is injective,
and the support of the coherent sheaf $\Cok(f)$ is $0$-dimensional.
By using the obstruction theory given in \cite{m1},
we can conclude that the filt scheme $\Filt(E,\vecl)$ is smooth,
and the tangent space is given by
$H^0(\Cok(f))$.
It is easy to see that
$\dim(H^0(\Cok(f)))=\rank(E)\cdot|\vecl|$.
In all, we have the following.
\begin{prop}
The variety $\Filt(E,\vecl)$ is 
smooth and of $\rank(E)\cdot|\vecl|$-dimension.
\hfill\qed
\end{prop}

The projection of of $\Filt(E,\vecl)\times C$
onto the $j$-th component by $q_j$.
We have the universal sequence of the quotients:
\[
 q_2^{\ast}(E)\lrarr \nbigf^u_h\lrarr \cdots
 \lrarr \nbigf^u_1.
\]
We put $\nbigk^u_j:=\ker(q_2^{\ast}(E)\lrarr \nbigf^u_j)$.
Then we obtain the following complex
over $\Filt(E,\vecl)\times C$:
\[
\begin{CD}
 \bigoplus_{i=1}^h
 Hom(\nbigk^u_i,\nbigk^u_{i})
@>{f^u}>>
 \bigoplus_{i=1}^h
 Hom(\nbigk^u_{i},\nbigk^u_{i-1}).
\end{CD}
\]
Then tangent bundle of $\Filt(E,\vecl)\times C$ is given by
the vector bundle $q_{1\,\ast}(\Cok(f^u))$.

\subsubsection{Torus action and the fixed point sets}
\label{subsubsection;1.9.85}

Let consider the case $E=\bigoplus_{\alpha=0}^{r-1}L_{\alpha}$,
for line bundles $L_{\alpha}$ $(\alpha=0,\ldots,r-1)$ over $C$.
Let $G_m^r$ denotes the $r$-dimensional torus.
Then we have the natural $G_m^r$-action on $E$,
which induces the $G_m^r$-action on $\Filt(E,\vecl)$.

For an $\vecl_{\ast}\in \Dec(\vecl,r)$,
we put as follows:
\[
 F^{\vecl}(\vecl_{\ast}):=
 \prod_{\alpha=0}^{r-1}
 Filt(L_{\alpha},\vecl_{\alpha})
\simeq
 \prod_{\alpha=0}^{r-1} C^{(\vecl_{\alpha})}.
\]
In particular, we use the notation
$F(\vecl_{\ast})$ and $F^Q(\vecl_{\ast})$
instead of
$F^{\vecc_l}(\vecl_{\ast})$ and $F^{l\cdot\vece_1}(\vecl_{\ast})$
for any $l$.

Let $\Filt(E,\vecl)^{G_m^r}$ denote
the fixed point set of $\Filt(E,\vecl_{\ast})$
with respect to the torus action above.
Then we have the natural isomorphism:
\[
 \Filt(E,\vecl)^{G_m^r}
\simeq
 \coprod_{\vecl_{\ast}\in\Dec(\vecl,r)}
 F^{\vecl}(\vecl_{\ast}).
\]

On $\Filt(E,\vecl)^{G_m^r}\times C$,
the universal sequence is decomposed as follows:
\[
 \bigoplus_{\alpha=0}^{r-1}
 \Bigl( 
 q_2^{\ast}L_{\alpha}
\lrarr
 \nbigf^u_{\alpha,h}
\lrarr
 \nbigf^u_{\alpha,h-1}
\lrarr\cdots\lrarr
 \nbigf^u_{\alpha,1}
 \Bigr).
\]
We put
$\nbigk^u_{\alpha,j}:=
 \ker(q_2^{\ast}L_{\alpha}\lrarr \nbigf^u_{\alpha,j})$
and $\nbigk^u_{\alpha,0}:=q_2^{\ast}L_{\alpha}$.
We put as follows:
\[
 \begin{array}{l}
 \nbign^{\vecl}_+(\vecl_{\ast}):=
\Cok\Bigl(
 \bigoplus_{\alpha<\beta}
 \bigoplus_{j=1}^h
 Hom(\nbigk^u_{\alpha,j},\nbigk^u_{\beta,j})
\lrarr
\bigoplus_{\alpha< \beta}
 \bigoplus_{j=1}^h
 Hom(\nbigk^u_{\alpha,j},\nbigk^u_{\beta,j-1})
 \Bigr),\\
\mbox{{}}\\
 \nbign^{\vecl}_-(\vecl_{\ast}):=
\Cok\Bigl(
 \bigoplus_{\alpha>\beta}
 \bigoplus_{j=1}^h
 Hom(\nbigk^u_{\alpha,j},\nbigk^u_{\beta,j})
\lrarr
\bigoplus_{\alpha>\beta}
 \bigoplus_{j=1}^h
 Hom(\nbigk^u_{\alpha,j},\nbigk^u_{\beta,j-1})
 \Bigr).
\end{array}
\]
Then we put
$\nbign^{\vecl}(\vecl_{\ast})
 =\nbign^{\vecl}_+(\vecl_{\ast})
 \bigoplus
  \nbign^{\vecl}_-(\vecl_{\ast})$.
Then $q_{2\,\ast}(\nbign^{\vecl}(\vecl_{\ast}))$
is the normal bundle of $F^{\vecl}(\vecl_{\ast})$.
We denote it by $N^{\vecl}(\vecl_{\ast})$.

\subsubsection{The decompositions due to Bianicki-Birula, Carrel and Sommes}

Let consider the sub-torus $G_m\subset G_m^r$
given by $t\longmapsto (t,t^2,\ldots,t^r)$.
Then we obtain the $G_m$-action
on $N^{\vecl}(\vecl_{\ast})$.
The negative part $N^{\vecl}_-(\vecl_{\ast})$ is given by
$q_{2\,\ast}(\nbign_-^{\vecl}(\vecl_{\ast}))$.
It is easy to see the equality
$\rank(N^{\vecl}_-(\vecl_{\ast}))=\co(\vecl_{\ast})$.

Due to the theory of Bianicki-Birula \cite{bb}
or Morse-theoretic consideration,
we have the decomposition of $\Filt(E,\vecl_{\ast})$
into the unstable manifolds:
\[
 Filt(E,\vecl_{\ast})=
 \coprod_{\vecl_{\ast}\in \Dec(\vecl,r)}
 F^{\vecl}(\vecl_{\ast})^+.
\]
It satisfies the following:
\begin{itemize}
\item
The complex codimension of $F^{\vecl}(\vecl_{\ast})^+$ is 
$\rank(N_-^{\vecl}(\vecl_{\ast}))=\co(\vecl_{\ast})$.
\item
We have the $G_m$-equivariant morphism
$\pi_{\vecl,\vecl_{\ast}}:
 F^{\vecl}(\vecl_{\ast})^+\lrarr F^{\vecl}(\vecl_{\ast})$.
It is an affine bundle.
\item
Let $x$ be a point of $F^{\vecl}(\vecl_{\ast})^+$.
Then we have
${\displaystyle
\lim_{G_m\ni t\to 0}t\cdot x=\pi_{\vecl,\vecl_{\ast}}(x)
 \in F^{\vecl}(\vecl_{\ast})}$.
\end{itemize}
When $\vecl=l\cdot\vece_1$,
the moduli theoretic meaning of $F^{l\cdot\vece_1}(\vecl_{\ast})^+$
is given in \cite{m1}.
Following Carrel and Sommes \cite{cs},
we consider the natural embedding
$F^{\vecl}(\vecl_{\ast})^+\lrarr F^{\vecl}(\vecl_{\ast})\times \Filt(E,\vecl)$.
The closure $\overline{F^{\vecl}(\vecl_{\ast})}^+$
induces the correspondence map:
\[
 \xi^{\vecl}(\vecl_{\ast};\cdot):
 H^{\ast}(F^{\vecl}(\vecl_{\ast}))[-2\co(\vecl_{\ast})]
\lrarr
 H^{\ast}(\Filt(E,\vecl)).
\]
Due to the result of Carrel and Sommes,
we have the decomposition:
\[
 H^{\ast}(\Filt(E,\vecl))
\simeq
 \bigoplus_{\vecl_{\ast}\in\Dec(\vecl,r)}
 H^{\ast}(F^{\vecl}(\vecl_{\ast}))[-2\co(\vecl_{\ast})].
\]

We have already given the bijection
$\Dec(\vecl,r)\simeq \nbigb(\vecl,r)$.
We identify them.
Then we have the morphism
$\xi^{\vecl}(\vecv_{\ast}:\cdot):
 H^{\ast}(F^{\vecl}(\vecv_{\ast}))[-2\co(\vecv_{\ast})]
\lrarr
 H^{\ast}(\Filt(E,\vecl))$,
and the decomposition:
\[
 H^{\ast}(\Filt(E,\vecl))
 =\bigoplus_{\vecl_{\ast}\in\nbigb(\vecl,r)}
 \Image(\xi^{\vecl}(\vecv_{\ast};\cdot))
 \simeq
 \bigoplus_{\vecv_{\ast}\in \nbigb(\vecl,r)}
 H^{\ast}(F^{\vecl}(\vecv_{\ast}))[-2\co(\vecv_{\ast})].
\]
By using the decomposition,
we introduce the following filtration:
\[
 F_hH^{\ast}(\Filt(E,\vecl))=
 \bigoplus_{
 \substack{
 \vecv_{\ast}\in \nbigb(\vecl,r)\\
 \co(\vecv_{\ast})\leq h}
 }
\Image(\xi^{\vecl}(\vecv_{\ast};\cdot))
\]
The filtration is depending on a decomposition
$E=\bigoplus_{\alpha=0}^{r-1}L_{\alpha}$.

%% file: 12.2.tex
\subsection{Complete filt scheme}

\label{subsection;1.9.80}

\subsubsection{An inductive description of complete filt schemes}

Let $q_j$ denote the projection of $\Filt(E,n)\times C$
onto the $j$-th component.
We have the universal sequence of quotients
$q_2^{\ast}(E)\lrarr \nbigf^u_n\lrarr \nbigf^u_{n-1}\lrarr\cdots\lrarr
 \nbigf^u_1$.
We put $\nbigk^u_n:=\ker(q_2^{\ast}(E)\lrarr\nbigf^u_n)$,
which is locally free on $\Filt(E,n)\times C$.
Thus we obtain the associated projective space bundle
$\proj(\nbigk_n^{u\,\lor})$.
We denote the associated tautological line bundle by $\nbigl_{n+1}$.
We see the following rather obvious lemma.
\begin{lem}
The variety 
$\Filt(E,n+1)$ is naturally isomorphic to
$\proj(\nbigk_n^{u\,\lor})$.
\end{lem}
\pf
Let $U$ be a scheme.
We denote the projection of $U\times C$ onto the $j$-th component
by $q_j$.
Let $q_2^{\ast}(E)
 \lrarr\nbigf_{n+1}\lrarr \nbigf_{n}\lrarr\cdots\lrarr\nbigf_{1}$
be a sequence of quotients such that
such that the length of $\ker(\nbigf_{j}\lrarr\nbigf_{j-1})$
is $1$.
We naturally obtain the sequence,
$q_2^{\ast}(E)
 \lrarr\nbigf_{n}\lrarr\cdots\lrarr\nbigf_{1}$.
It gives the morphism $g:U\lrarr \Filt(E,n)$
such that $(g\times id_C)^{\ast}\nbigf_j^u=\nbigf_j$.
Here $g\times id_C$ is naturally induced morphism
$U\times C\lrarr \Filt(E,n)\times C$.
We put as follows:
\[
 \nbigk_n:=\ker(q_2^{\ast}(E)\lrarr\nbigf_n)=
 (g\times id_C)^{\ast}\nbigk_n^u,
\quad
\quad
 \nbigg_{n+1}:=\ker(\nbigf_{n+1}\lrarr\nbigf_n).
\]
Then we have the naturally induced morphism
$\nbigk_n\lrarr\nbigg_{n+1}$ defined over $U\times C$.
It induces the section
$U\lrarr \proj(\nbigk_n^{\lor})$.
As a result, we obtain the morphism
$U\lrarr \proj(\nbigk_n^{u\,\lor})$.
In particular,
we obtain the morphism
$\Phi_{n+1}:\Filt(E,n+1)\lrarr\proj(\nbigk_n^{u\,\lor})$.

Let $\pi_n$ denote the projection of $\proj(\nbigk_n^{u\,\lor})\lrarr
 \Filt(E,n)\times C$.
We put $\chi_n:=q_1\circ \pi_n$ and $\rho_n:=q_2\circ \pi_n$.
Then we obtain the morphism
$\chi_n\times id_C:\projnbigku\times C\lrarr \Filt(E,n)\times C$.
Then we have the following morphisms:
\[
 (\chi_n\times id_C)^{\ast}\nbigk_n^u
\lrarr
 q_2^{\ast}E
\lrarr
 (\chi_n\times id_C)^{\ast}\nbigf_n^u
\lrarr\cdots\lrarr
 (\chi_n\times id_C)^{\ast}\nbigf_1^u.
\]

On $\projnbigku$, we have the naturally defined morphism
$\psi_n:\chi_n^{\ast}\nbigk^u_n\lrarr\nbigl_{n+1}$.
We have the closed embedding
$(id\times \rho_n):\projnbigku\lrarr\projnbigku\times C$.
Let $X_n$ denote the image.
Then we obtain the morphism
$\lambda_n:
(\chi_n\times id_C)^{\ast}\nbigk^u_n\otimes\nbigo_{X_n}
\lrarr
 q_1^{\ast}\nbigl_{n+1}\otimes\nbigo_{X_n}$,
which is induced by $\psi_n$ and the isomorphism $X_n\simeq \projnbigku$.
We put $\nbigk^{u}_{n+1}:=\ker(\lambda_n)$,
and $\nbigf^u_{n+1}:=\cok(\nbigk_{n+1}^u\lrarr q_2^{\ast}E)$.
We denote
$(\chi_{n+1}\times id)^{\ast}\nbigf_{j}^u$
by $\nbigf_j^u$ for $j=1,\ldots,n$.
Then we obtain the sequence
$q_2^{\ast}(E)\lrarr\nbigf^u_{n+1}\lrarr\nbigf^u_{n}\lrarr\cdots
\lrarr \nbigf^u_1$
defined over $\projnbigku\times C$.
In particular, we obtain the morphism
$\Psi_{n+1}:\projnbigku\lrarr \Filt(E,n+1)$.

It is easy to see that $\Phi_{n+1}$ and $\Psi_{n+1}$
are inverses each other.
Thus $\Filt(E,n+1)$ and $\projnbigku$ are isomorphic
via the morphisms.
\hfill\qed

\vspace{.1in}
We will not distinguish them in the following.
By our construction,
we also know the following:
\begin{cor}
The sheaf $\nbigg^u_{n+1}:=\ker(\nbigf^u_{n+1}\lrarr\nbigf^{u}_{n})$
is isomorphic to
$q_1^{\ast}\nbigl_{n+1}\otimes \nbigo_{X_n}$.
In particular, we have
$q_{1\,\ast}\bigl(\nbigg^u_{n+1}\bigr)\simeq \nbigl_{n+1}$.
\hfill\qed
\end{cor}

\subsubsection{The structure of $H^{\ast}(Filt(E,n),\seisuu)$}

In the rest of the subsection \ref{section;1.13.20},
we consider the cohomology ring with $\seisuu$-coefficient.
\label{subsubsection;1.9.81}
For any $m\leq n$,
we have the natural morphism
$\eta_{n,m}:\Filt(E,n)\lrarr \Filt(E,m)$
defined by the following functor:
\[
 \Bigl(E\rarr \nbigf_n\rarr\cdots
  \rarr\nbigf_m\rarr\cdots
 \rarr \nbigf_1
 \Bigr)
\longmapsto
 \Bigl(E\rarr \nbigf_m\rarr\cdots\rarr \nbigf_1
 \Bigr).
\]
We clearly have
$\eta_{n,m}:=
 \chi_{m}\circ \chi_{m+1}\circ\cdots \circ \chi_{n-1}$.
We denote $\eta_{n,m}^{\ast}\nbigl_m$ by $\nbigl_m$ for simplicity
of notation.

For the universal filtration
$q_2^{\ast}(E)\lrarr\nbigf^u_n\lrarr\cdots\lrarr\nbigf^u_{1}$,
we put $\nbigg^u_j:=\ker(\nbigf^u_{j}\lrarr\nbigf^u_{j-1})$.
Then the support of $\nbigg^u_j$ induces the morphism
$p_j:\Filt(E,n)\lrarr C$.
Note that $p_{n}$ is same as $\rho_{n-1}$ in the previous subsection.
In all, we have the morphism
$p=\prod_{j=1}^np_j:\Filt(E,n)\lrarr C^n$.
We have the induced morphism
$p^{\ast}:H^{\ast}(C^n)\lrarr H^{\ast}(\Filt(E,n))$.
For any element $a\in H^{\ast}(C^n)$,
the element $p^{\ast}(a)$ is denoted by $a$ for simplicity of notation.

We have the morphism
$\rho_{n-1}\times id_C:\Filt(E,n)\times C\lrarr C\times C$.
Then $X_{n-1}$ is same as $(\rho_{n-1}\times id_C)^{\ast}\Delta$.
It is reworded as follows:
We have the morphism
$p\times id_C:\Filt(E,n)\times C\lrarr C^n\times C$.
Then $X_{n-1}$ is same as $(p\times id_C)^{\ast}\Delta_{n,n+1}$.
In particular,
the sheaf $\nbigg^u_{n}$ over $\Filt(E,n)\times C$
is isomorphic to the following:
\[
 q_1^{\ast}\nbigl_n\otimes
\Cok\Bigl(
 \nbigo(-\Delta_{n,n+1})
\lrarr
 \nbigo
 \Bigr).
 \]
Here we denote the pull back of the divisor $\Delta_{n,n+1}$
via the morphism $p\times id_C$ by the same notation $\Delta_{n,n+1}$.
Hence 
the sheaf $\nbigg^u_{j}$ over $\Filt(E,n)\times C$
is isomorphic to the following:
\[
 q_1^{\ast}\nbigl_j\otimes
\Cok\Bigl(
 \nbigo(-\Delta_{j,n+1})
\lrarr
 \nbigo
 \Bigr).
\]
In the $K$-theory of the coherent sheaves
on $\Filt(E,n)$,
we have the equality
$ \nbigku_n=
 q_2^{\ast}(E)-\sum_{j=1}^n\nbiggu_j$.
Hence we obtain the following equality for the total Chern classes:
\[
  c_{\ast}\bigl(\nbigku_n\bigr)=
 c_{\ast}\bigl(
 q_2^{\ast}(E)\bigr)
\times
 \prod_{j=1}^n
 c_{\ast}\bigl(
 \nbiggu_j
 \bigr)^{-1}.
\]

We denote the first Chern class of $\nbigl_j$ by $\omega_j$.
Note that the composition
$\Filt(E,n+1)\lrarr \Filt(E,n)\times C\lrarr C$
is same as $p_{n+1}:\Filt(E,n+1)\lrarr C$.
Thus we obtain the following relation in the cohomology ring
$H^{\ast}(\Filt(E,n+1))$:
\[
 \sum_{j=0}^r a_j\cdot \omega_{n+1}^{r-j}=0,
\quad
 \sum_{j=0}^r a_j=
 p_{n+1}^{\ast}\bigl(
 c_{\ast}(E)\bigr)\cdot
 \prod_{j=1}^n\frac{1+\omega_j-\Delta_{j\,n+1}}{1+\omega_j}.
\]
Due to the general theory for the projective space bundles,
the relation determines the structure of $H^{\ast}(\Filt(E,n+1))$
over $H^{\ast}(\Filt(E,n))$ with the generator $\omega_{n+1}$.

By an inductive argument, we obtain the following proposition.
\begin{prop}
The ring $H^{\ast}(\Filt(E,n))$ is the quotient ring of
$H^{\ast}(C^n)[\omega_1,\ldots,\omega_n]$ divided by the following
relations for $1\leq h\leq n$:
\[
 \sum_{j=0}^r a_j^{(h)}\cdot \omega_h^{r-j}.
\]
Here $a_j^{(h)}$ is described
in terms of $\Delta_{i,h}$ and $\omega_j$ $(i,j\leq h-1)$
as follows:
\[
 \sum_{j=0}^r a_j^{(h)}=
\prod_{i=1}^{h-1}\frac{1+\omega_i-\Delta_{i\,h}}{1+\omega_i}.
\]
\hfill\qed
\end{prop}

\subsubsection{Some easy corollaries}

As a corollary derived from the description above,
we know the structure of
the cohomology ring of infinite complete filt schemes.
Let $L_{\alpha}$ $(\alpha=0,1,2\ldots)$ be line bundles over $C$.
We put $E_{\beta}:=\bigoplus_{\alpha\leq \beta}L_{\alpha}$,
and $E_{\infty}:=\bigoplus_{\alpha}L_{\alpha}$.
Then we have the naturally defined inclusions
$\Filt(E_{\beta},n)\lrarr \Filt(E_{\beta+1},n)$.
Then we can consider the limit $\Filt(E_{\infty},n)$.

\begin{cor}
$H^{\ast}(\Filt(E_{\infty},n))$
is the polynomial ring over $H^{\ast}(C^n)$
generated by the first Chern classes $\omega_i$ for $i=1,\ldots,n$
\hfill\qed
\end{cor}

In general,
we have the classifying maps
$\mu_i:\Filt(E,n)\lrarr \proj^{\infty}$
corresponding to the line bundles $\nbigl_i\,\,(i=1,\ldots,n)$,
for any locally free sheaf $E$ over $C$.
In all we have the morphism
$\mu=\prod_{i=1}^n\mu_i: \Filt(E,n)\lrarr (\proj^{\infty})^n$.
Thus we obtain the morphism
$p\times \mu:\Filt(E,n)\lrarr C^n\times (\proj^{\infty})^n$.

\begin{cor}
Let consider the case $C=\proj^1$.
Then the morphism $\Filt(E_{\infty},n)\lrarr
 \bigl(\proj^1\times \proj^{\infty}\bigr)^n$
is homotopy equivalent.
\hfill\qed
\end{cor}

Let $E$ be a locally free coherent sheaf over $C$.
Recall that
the sheaf $\nbigg^u_i$ defined
over $\Filt(E,n)\times C$ is isomorphic to
$q_1^{\ast}\nbigl_i\otimes
 \Cok\Bigl(
 \nbigo(-\Delta_{i,n+1})\rarr \nbigo
  \Bigr)$.
Here $\Delta_{i,n+1}$ denotes the pull back of the divisor
$\Delta_{i,n+1}$ of $C^n\times C$.
Hence the total Chern class of $\nbigg_i$ is described as follows:
\begin{equation}
 c_{\ast}(\nbigg_i)=
 \frac{1+\omega_i}{1+\omega_i-\Delta_{i,n+1}}=
 \sum_{j=0}^{\infty}\left(\frac{\Delta_{i,n+1}}{1+\omega_i}\right)^j
=
1+\Delta_{i,n+1}\sum_{h=0}^{\infty}(-\omega_i)^h
 +\Delta_{i,n+1}^2(1+\omega_i)^{-2}.
\end{equation}
In particular,
we have $c_1(\nbigg_i)=\Delta_{i,n+1}$
and
$c_2(\nbigg_i)=-\Delta_{i,n+1}\cdot \omega_i+(2-2g(C))\cdot \pt_{i,l+1}$.
Here $g(C)$ denotes the genus of $C$.

The $j$-th Chern class $c_j(\nbigg_i)$
induces the correspondence map
$\psi_j^{(i)}:
 H^{\ast}(C)[-2(j-1)]\lrarr H^{\ast}(\Filt(E,n))$.
The following equalities are easily checked:
\begin{equation} \label{eq;1.9.1}
 \psi_1^{(i)}(a)=p_i^{\ast}(a),
\quad\quad
 \psi_2^{(i)}(1)=-\omega_i+(2-2g(C))\cdot p_i^{\ast}(pt).
\end{equation}
\begin{prop}
Thus the following set generates the ring $H^{\ast}(Filt^{comp}(E,l))$
over $\rnum$:
\[
 \bigcup_{i=1}^n
\left(\Image \psi_1^{(i)}
 \cup\{\psi_2^{(i)}(1)\}
\right).
\]
\end{prop}
\pf
The claim follows from (\ref{eq;1.9.1}).
\hfill\qed

The proposition will be generalized in the case of the other filt schemes.

%% file: 12.3.tex
\section{Some products in
 $H^{\ast}(\Filt(\bigoplus_{\alpha=0}^{r-1}L_{\alpha},n))$}

\subsection{Preliminary}
\label{subsection;1.14.1}

We have already known the ring structure of
$H^{\ast}\big(\Filt(\bigoplus_{\alpha=0}^{r-1}L_{\alpha},n)\big)$,
which is described by $\omega_i$ $(i=1\ldots,n)$.
On the other hand, we have the naturally defined cohomological classes
$\big\{\xi^{\vecc_n}(\vecv_{\ast};a)\,
 \big|\,\vecv_{\ast}\in \nbigb(\vecc_n,r),\,
 a\in H^{\ast}(F^{\vecc_n}(\vecv_{\ast}))\big\}$.
We calculate the action of $\omega_n$ on
$\xi^{\vecc_n}(\vecv_{\ast};a)$.
In the following, we put $E=\bigoplus_{\alpha=0}^{r-1}L_{\alpha}$.

\subsubsection{Replacement of notation}

Recall that $\nbigb(\vecc_n,r)$ is naturally isomorphic
to $[0,r-1]^n$.
Thus we use the notation $\vecv=(v_1,\ldots,v_n)$
instead of $\vecv_{\ast}=(\vecv_1,\ldots,\vecv_n)$.
Moreover we use $F(\vecv)$ instead of $F^{\vecc_n}(\vecv)$.
Similarly
we use $\xi(\vecv)$ instead of $\xi^{\vecc_n}(\vecv;1)$.
Note that we have $\xi^{\vecc_n}(\vecv;a)=a\cdot \xi(\vecv)$.
Then we have the decomposition:
\[
 H^{\ast}(\Filt(E,n))=
 \bigoplus _{\vecv\in [0,r-1]^n}
 H^{\ast}(C^n)\cdot \xi(\vecv).
\]

We will use the equivariant cohomology group
$H^{\ast}_{G_m^r}(\Filt(E,n))$.
Since the action of $G_m^r$ preserves the subvarieties
$F(\vecv)^+$, the class $\xi(\vecv)$ is naturally  lifted to
the equivariant cohomology class.
We denote them by the same notation.
Moreover, they give
the base of $H^{\ast}_{G_m^r}(\Filt(E,n))$
over the ring $H^{\ast}(C^n)[t_0,\ldots,t_{r-1}]$.

We also have the equivariant first Chern class
of $\nbigl_i$,
which is denoted also by $\omega_i$.

The inclusion $F(\vecv)\subset\Filt(E,n)$ induces
the naturally defined ring morphism:
\[
 \res(\vecv):
 H^{\ast}_{G_m^r}(\Filt(E,n))\lrarr H^{\ast}_{G_m^r}(F(\vecv))
\simeq H^{\ast}(F(\vecv))[t_0,\ldots,t_{r-1}].
\]
For an element $f$ of $H^{\ast}_{G_m^r}(\Filt(E,n))$,
the element $\res(\vecv)(f)\in H^{\ast}(F(\vecv))[t_0,\ldots,t_{r-1}]$
is denoted by $f^{\vecv}$.

\subsection{Some lemmas}

We prepare some lemmas which are obtained by some geometric considerations.

\begin{lem}
We have the following equality:
\[
 \omega_i^{\vecv}=
 t_{v_i}-\sum_{\substack{k<i\\v_k=v_i}} \Delta_{k,i}
 +p_i^{\ast}c_1(L_{v_i}).
\]
\end{lem}
\pf
We only have to show the equality in the case $i=n$.
Recall that $\Filt(E,n)$ is
isomorphic to the associated projective space bundle
with $\nbigk_{n-1}^u$ over $\Filt(E,n-1)\times C$,
and $\nbigl_n$ is the tautological line bundle under the identification.
Then the restriction $\nbigl_{n\,|\,F(\vecv)}$
is isomorphic to
$\ker(q_2^{\ast}L_{v_n}\lrarr \nbigf_{v_n,n-1})$.
It is isomorphic to the following line bundle:
\[
 q_2^{\ast}L_{v_n}\otimes
 \nbigo\Bigl(
 -\sum_{i<n,\,v_i=v_n}
 \Delta_{i\,n}
 \Bigr).
\]
Then we obtain the result.
Note the composition
$\Filt(E,n)\lrarr \Filt(E,n-1)\times C\stackrel{q_2}{\lrarr} C$
is same as $p_n$.
\hfill\qed

\begin{df}
We introduce the order $\leqzero$ on $[0,r-1]^n$
defined as follows:
\begin{quote}
 Let $\vecv=(v_1,\ldots,v_n)$ and $\vecv'=(v_1',\ldots,v_n')$
 be elements of $[0,r-1]^n$.
Then $\vecv\leqzero \vecv'$ if and only if
$v_i\leq v_i'$ for any $i=1,\ldots,n$.
\end{quote}
By the inclusion $\nbigb(n,r)\subset [0,r-1]^n$,
we have the induced order on $\nbigb(n,r)$.
We denote the order by the same notation.
\hfill\qed
\end{df}

Recall that $[f]$ denotes the top term
for $f\in H^{\ast}(F(\vecv))[t_0,\ldots,t_{r-1}]$
(the subsection \ref{subsection;1.13.30}).
\begin{lem} \label{lem;1.9.20}
Assume that  $\vecw\in[0,r-1]^n$ satisfies 
$\vecw\geq_0\vecv$.
Then we have the following equality:
\[
 [\xi(\vecv)^{\vecw}]=
 \prod_{j=1}^n\prod_{i=1}^{v_j-1}(t_{w_j}-t_i).
\]
In particular, we have $\deg_t(\xi(\vecv)^{\vecw})=\co(\vecv)$.
Reversely, $\deg_t(\xi(\vecv)^{\vecw})=\co(\vecv)$
implies $\vecw\geq_0\vecv$.
\end{lem}
\pf
Let $P=(x_1,\ldots,x_n)$ be a point of $C^n$
such that $x_i\neq x_j$ $(i\neq j)$.
Then we have the closed embedding
$\Filt(E,n)\times_{C^n}\{P\}\lrarr \Filt(E,n)$.
Note that $\Filt(E,n)\times_{C^n}\{P\} $ is isomorphic to
$(\proj^{r-1})^n$.
We only have to consider the image of $\xi(\vecv)$
to $H^{\ast}(\Filt(E,n)\times_{C^n}\{P\})$.
\hfill\qed

\begin{lem}\label{lem;1.15.10}
If $\xi(\vecv)^{\vecw}\neq 0$,
then we have an element $\sigma\in \gbigs_n$
such that $\sigma(\vecv)\leq \vecw$.
\end{lem}
\pf
Let consider the morphism $\Psi:\Filt(E,n)\lrarr \Quot(E,n)$.
Then $\Psi(F(\vecv))$ is contained in $F^Q(\Nor(\vecv))$,
and $\Psi(\overline{F(\vecv)}^+)\subset \overline{F^Q(\Nor(\vecv))}^+$.
In \cite{m1}, we obtained that
$\overline{F^Q(\vecu)}^+ \cap F^Q(\vecu')\neq \emptyset
\Longleftrightarrow \vecu\leqzero \vecu'$.
Thus we are done.
\hfill\qed

\begin{lem} \label{lem;1.9.2}
Let $\vecw$ be an element of $[0,r-1]^n$,
such that there exists $\sigma\in\gbigs_n$ satisfying
$\sigma(\vecv)=\vecw$.
We put $I:=\{i\,|\,v_i\neq w_i\}$.
Then we have the following inequality:
\[
 \deg_t\xi(\vecv)^{\vecw}\leq co(\vecv)-|I|+1.
\]
\end{lem}
\pf
Let consider the restriction of $\xi(\vecv)$
to $\Filt(E,n)\times_{C^n}(C^n-\Delta_I)$.
It is easy to see the following:
\[
 \overline{F(\vecv)^+}\cap F(\vecw)\cap
 \bigl[
  \Filt(E,n)\times_{C^n}(C^n-\Delta_I)
 \bigr]=\emptyset.
\]
Thus we are done.
\hfill\qed

Let $\vecw$ be an element of $[0,r-1]^n$.
We put
$J(\vecv,\vecw):=\{\sigma\in\gbigs_n\,|\,\sigma(\vecv)\leq \vecw\}$.
We denote the projection of $[0,r-1]^n$ onto the $i$-th component
by $q_i$.
For any element $\sigma\in J(\vecv,\vecw)$,
we put $I(\sigma):=\{i\,|\,q_i(\vecv)\neq q_i(\sigma(\vecv))\}$.
Then we obtain the number
$d(\vecv,\vecw):=
 \min\bigl\{|I(\sigma)|\,\big|\,\sigma\in J(\vecv,\vecw)\bigr\}$.

\begin{lem}
We have the following inequality:
\[
 \deg_t(\xi(\vecv)^{\vecw})
\leq
 \co(\vecv)-d(\vecv,\vecw)+1.
\]
\end{lem}
\pf
Similar to Lemma \ref{lem;1.9.2}.
\hfill\qed

\begin{cor} \label{cor;1.16.1}
Let $\vecw$ be an element of $[0,1]$ satisfying
$\co(\vecw)=\co(\vecv)$ and
$\deg_t(\xi(\vecv)^{\vecw})=co(\vecv)-1$.
Then there exists a transposition $\tau$
such that $\tau(\vecv)=\vecw$.
Namely, there exist $i\neq j$ satisfying
$v_i=w_j$, $v_j=w_i$ and $v_k=w_k$ $(k\neq i,j)$.
\hfill\qed
\end{cor}

\subsection{Calculation}


We consider the case $E=\bigoplus_{\alpha=0}^{r-1}L_{\alpha}$
for line bundles $L_{\alpha}$ on $C$.

\begin{thm} \label{thm;1.9.40}
Let $\vecv$ be an element of $[0,r-1]^n$.
In the equivariant cohomology ring
$H^{\ast}_{G_m^r}(\Filt(E,n))$,
we have the following equality:
\begin{equation} \label{eq;1.9.15}
\xi(\vecv+\vece_n)=
 (\omega_n-p_n^{\ast}L_{v_n}-t_{v_n})\cdot \xi(\vecv)
+\sum_{\substack{k<n,\\ v_k\leq v_n}}
 \Delta_{k,n}\cdot \xi(\tau_{k,n}(\vecv)).
\end{equation}
As a direct corollary,
we have the following equality in the non-equivariant cohomology ring
$H^{\ast}(\Filt(E,n))$:
\begin{equation} \label{eq;1.9.71}
 \xi(\vecv+\vece_n)=
 (\omega_n-p_n^{\ast}L_{v_n})\cdot \xi(\vecv)
+\sum_{\substack{k<n,\\ v_k\leq v_n}}
 \Delta_{k,n}\cdot \xi(\tau_{k,n}(\vecv)).
\end{equation}
Here $\tau_{k\,n}$ denotes the transposition of $k$ and $n$.
\end{thm}
\pf
We put $X:=(\omega_n-\omega_n^{\vecv})\cdot \xi(\vecv)$.
Then we only have to prove the following equality:
\[
 X=
 \sum_{\substack{k<n,\\ v_k<v_n}} -\Delta_{k,n}\cdot \xi(\tau_{k,n}\vecv)
 +\xi(\vecv+\vece_n).
\]
Note $X^{\vecv}\neq 0$. First we see the following lemma.
\begin{lem} \label{lem;1.9.10}
The element $X$ has the following description:
\[
 X=
\sum_{\substack{i<n\\ v_i\neq v_n}}A_i\cdot \xi(\tau_{i\,n}(\vecv))
+\xi(\vecv+\vece_n).
\]
Here $A_i$ are elements of $H^2(C^n)$.
\end{lem}
\pf
We put 
$\nbigd_1:=\{\vecw\in \nbigb(n,r)
 \,|\,\exists \sigma\in\gbigs_n,\,X^{\sigma(\vecw)}\neq 0
 \}$.
Then we know that $\vecw\geq_0 \Nor(\vecv)$ for any $\vecw\in \nbigd_1$
(Lemma \ref{lem;1.15.10}).

We have the description $X=\sum_{\vecu\in[0,r-1]^n} a(\vecu)\cdot \xi(\vecu)$.
We put
$\nbigd_2:=\{\vecu\in \nbigb(n,r)
 \,|\, \exists \sigma\in\gbigs_n,\,
 a(\sigma(\vecu))\neq 0
 \}$.
Let $\vecw_0$ be a minimal element of $\nbigd_2$
with respect to the order $\leqzero$.
Let $\vecw_0'$ be a reordering of $\vecw_0$
satisfying the following:
\begin{equation} \label{eq;1.9.5}
 \deg_t a(\sigma\vecw_0')\leq
 \deg_t a(\vecw_0'),
\quad \mbox{for any } \sigma\in\gbigs_n.
\end{equation}
Since $\vecw_0$ is minimal, we have the following equality:
\[
 X^{\vecw_0'}=
 \sum_{\sigma\in \gbigs_n} a(\sigma\vecw_0')\cdot
 \xi(\sigma\vecw_0')^{\vecw_0'}.
\]
Due to the inequalities (\ref{eq;1.9.5})
and Lemma \ref{lem;1.9.2},
we know that the right hand side is not $0$.
Thus $\vecw_0$ is contained in $\nbigd_1$.
In particular, $\vecw_0\geq \Nor(\vecv)$.

Thus $X$ has the following description:
\[
 X=\sum_{
 \substack{\Nor(\vecu)\geq \Nor(\vecv)}}
 a(\vecu)\cdot \xi(\vecu).
\]
Let $\vecv'$ be a reordering of $\vecv$,
satisfying the following:
\begin{equation} \label{eq;1.9.6}
 \deg_t a(\sigma(\vecv'))
\leq
 \deg_t a(\vecv')
\quad
 \mbox{ for any } \sigma\in\gbigs_n.
\end{equation}
Then we obtain the following equality,
by using the inequality
$\deg_t\xi(\sigma\vecv')^{\vecv'}<\deg_t \xi(\vecv')^{\vecv'}$
$(\sigma\neq 1)$:
\[
 \deg_t X^{\vecv'}=
 \deg_t a(\vecv')+\deg_t \xi(\vecv')^{\vecv'}
=\deg_t a(\vecv')+\co(\vecv).
\]
Assume that $\vecv'=\vecv$, and then 
it contradicts $X^{\vecv}=0$.
Thus we may assume that $\vecv'\neq\vecv$.
Then we obtain the following inequalities:
\[
 co(\vecv)\geq
 1+\deg_t\xi(\vecv)^{\vecv'}
\geq
 \deg_t(\omega_n^{\vecv'}-\omega_n^{\vecv})+\deg_t\xi(\vecv)^{\vecv'}
=\deg_t X^{\vecv'}=\co(\vecv)+\deg_t a(\vecv').
\]
It implies the following:
\begin{enumerate}
\item \label{number;1.9.7}
 $\deg_t a(\vecv')=0$.
 By seeing the cohomological degree,
 we can conclude that $a(\vecv')\in H^2(F(\vecv'))$.
\item \label{number;1.9.8}
 $\deg_t(\xi(\vecv)^{\vecv'})=-1+\co(\vecv)$.
Due to Corollary \ref{cor;1.16.1},
we have the transposition $\tau$ such that $\tau(\vecv')=\vecv$. 
\item \label{number;1.9.9}
 $\deg_t(\omega_n^{\vecv'}-\omega_n^{\vecv})=1$.
If $\tau(n)=n$,
then we have $\deg_t(\omega_n^{\vecv'}-\omega_n^{\vecv})=0$.
Hence we know that $\tau(n)\neq n$.
\end{enumerate}
Due to \ref{number;1.9.7},
the assertions \ref{number;1.9.8} and \ref{number;1.9.9}
hold for any reordering of $\vecv''$
such that $a(\vecv'')\neq 0$.

For any reordering $\vecv'$ of $\vecv$ such that $\vecv'\neq\vecv$,
we have $\deg_t a(\vecv')=0$.
If $a(\vecv)\neq 0$, then we obtain
$\deg_t X^{\vecv}\geq \co(\vecv)$,
which contradicts $X^{\vecv}=0$.
Thus we obtain $a(\vecv)=0$.

In all, $X$ has the following description:
\[
 X=\sum_{\substack{i<n\\ v_i\neq v_n}}A_i\cdot \xi(\tau_{i\,n}\vecv)
+\sum_{N(\vecu)>_0N(\vecv)}
 a(\vecu)\cdot \xi(\vecu).
\]
We put $Y=X-\sum_{i<n}A_i\cdot \xi(\tau_{i\,n}\vecv)=
 \sum_{\vecu}b(\vecu)\cdot \xi(\vecu)$.

We put
$\nbigd_3:=
 \{\vecw\in\nbigb(n,r)\,|\,
    \exists \sigma\in\gbigs_n,\,b(\sigma\vecw)\neq 0\}$.
We already have $\Nor(\vecv)<_0\vecw$ for any $\vecw\in \nbigd_4$.
In particular, $\co(\vecw)\geq \co(\vecv)+1$.
On the other hand, we know  $\co(\vecw)\leq \co(\vecv)+1$
due to the cohomological degree of $Y$.
Thus we know $\co(\vecw)=\co(\vecv)+1$.

Then we obtain that the cohomological degree
of $a(\vecw)$ is at most $0$.
To know such a value,
we only have to compare the both sides
in the cohomology ring $H^{\ast}(\Filt(E,n)\times _{C^n}\{P\})$
for any point $P=(x_1,\ldots,x_n)$ such that
$x_i\neq x_j$ $(i\neq j)$.

Hence the proof of Lemma \ref{lem;1.9.10} is completed.
\hfill\qed

\vspace{.1in}
Let us return to the proof of Theorem \ref{thm;1.9.40}.
We use an induction.
The following claim is called $Q(\gamma)$.
\begin{quote}
$Q(\gamma)$:
If $v_n\leq \gamma$,
then the equality (\ref{eq;1.9.15}) holds.
\end{quote}

The following claim is called $P(\gamma)$.
\begin{quote}
$P(\gamma)$:
 If $v_n\leq \gamma$, the following holds:
\begin{itemize}
 \item If $v_n>v_k$, we have the following equality for any $m\geq 0$:
 \begin{equation}
 \bigl[
 \xi(\vecv)^{\tau_{n,k}\vecv+m\cdot \vece_k}
 \bigr]
=\Delta_{k\,n}\cdot \prod_{j\neq n}\prod_{0\leq i\leq v_j-1}
 (t_{v_j}-t_{i})
\times
 \prod_{\substack{0\leq i\leq v_n-1\\i\neq v_k}}
 (t_{v_n+m}-t_i).
 \end{equation}
\item
If $v_n=v_k$, then we have the following equality for any $m\geq 0$:
\begin{equation}\label{eq;1.9.31}
 \bigl[
 \xi(\vecv)^{\vecv+m\cdot \vece_k}\bigr]
=
 \prod_{j\neq k}\prod_{i=0}^{v_j-1}(t_{v_j}-t_i)
\times
 \prod_{i=0}^{v_k-1}(t_{v_k+m}-t_i).
\end{equation}
\item
 If $v_n<v_k$ and if $0\leq m< v_k-v_n$,
then $\deg_t\xi(\vecv)^{\tau_{n,k}(\vecv)+m\cdot \vece_k}
 \leq \co(\vecv)-2$.
\end{itemize}
\end{quote}

\begin{lem}
The claim $P(0)$ holds.
\end{lem}
\pf
When $v_k=v_n=0$, the claim follows from Lemma \ref{lem;1.9.20}.
Let consider the case $v_k>v_n=0$.
The class $\xi(\vecv)$ is defined over $\Filt(E,n-1)$.
Namely it is the pull back of the corresponding class
via the morphism $\eta_{n,n-1}:\Filt(E,n)\lrarr \Filt(E,n-1)$.
Thus, $\xi(\vecv)^{\vecu+l\cdot \vece_n}=0$ for some $l\neq 0$
if and only if $\xi(\vecv)^{\vecu}=0$.
We know $\xi(\vecv)^{\tau_{n,k}\vecv-v_k\vece_n}=0$,
since $\Nor(\vecv)\not\leq_0\Nor(\vecv-v_k\cdot \vece_n)$.
Thus we obtain the vanishing $\xi(\vecv)^{\tau_{n,k}\vecv}=0$.
\hfill\qed

\begin{lem}
$Q(\gamma-1)+P(\gamma)\Longrightarrow Q(\gamma)$.
\end{lem}
\pf
We only have to show the equality (\ref{eq;1.9.15})
when $v_n=\gamma$.
We already have the equality:
\[
 (\omega_n-\omega_n^{\vecv})\cdot  \xi(\vecv)
=\xi(\vecv+\vece_n)+
\sum_{\substack{i<n\\v_i\neq v_n}} A_i \cdot\xi(\tau_{i,n}\vecv).
\]
If $k\neq i$, then we have the inequality
$\deg_t\xi(\tau_{k\,n}\vecv)^{\tau_{i\,n}\vecv}<\co(\vecv)-1$.
Thus we have the following equalities for any $k<n$
such that $A_k\neq 0$:
\begin{equation} \label{eq;1.9.21}
  \bigl[
 \omega_n^{\tau_{k,n}\vecv}-\omega_n^{\vecv}\bigr]
 \cdot
 \bigl[
 \xi(\vecv)^{\tau_{n,k}\vecv}
 \bigr]
=A_k\cdot [\xi(\tau_{k,n}\vecv)^{\tau_{k,n}\vecv}].
\end{equation}

When $v_k\neq v_n$,
we have the equality
$[(\omega_n-\omega_n^{\vecv})^{\tau_{k\,n}\vecv}]=
 t_{v_k}-t_{v_n}$.
We have the following equality,
due to Lemma \ref{lem;1.9.20}:
\[
 \big[\xi(\tau_{n,k}\vecv)^{\tau_{n,k}\vecv}\big]
=\prod_{j=1}^n\prod_{i=0}^{v_j-1}(t_{v_j}-t_i).
\]
On the other hand,
if $v_k>v_n$, then $P(\gamma)$
implies that $\deg_t\xi(\vecv)^{\tau_{n,k}\vecv}<\co(\vecv)-1$.
Thus the equation (\ref{eq;1.9.21}) cannot hold in this case.
Namely $A_k=0$ if $v_k>v_n$.

Let consider the case $v_k<v_n$.
Then $P(\gamma)$ implies the following equality:
\[
 [\xi(\vecv)^{\tau_{n,k}\vecv}]
=\Delta_{n,k}\cdot
 \prod_{j=1}^{n-1}\prod_{i=0}^{v_j-1}(t_{v_j}-t_i)
\times
 \prod_{\substack{0\leq i\leq v_{n-1}\\i\neq v_k}}
 (t_{v_n-t_i}).
\]
Then we obtain the equality $A_k=-\Delta_{n,k}$ in this case.
Thus we obtain $Q(\gamma)$.
\hfill\qed

\begin{lem} \label{lem;1.9.35}
$Q(\gamma)+P(\gamma)\Longrightarrow P(\gamma+1)$.
\end{lem}
\pf
Due to the hypothesis $Q(\gamma)$, we have the equality:
\begin{multline}\label{eq;1.9.30}
 \xi(\vecv+\vece_n)^{\tau_{n,k}(\vecv+\vece_n)+m\cdot \vece_k}
=
 \bigl(
 \omega_n^{\tau_{n,k}\vecv+(m+1)\cdot\vece_k}
-\omega_n^{\vecv}
 \bigr)\cdot
 \xi(\vecv)^{\tau_{n,k}\vecv+(m+1)\cdot\vece_k}\\
+\sum_{j<n,v_j<v_n}
 \Delta_{n,j}\cdot
  \xi(\tau_{n,j}\vecv)^{\tau_{n,k}\vecv+(m+1)\vece_k}.
\end{multline}

We divide the claim into the following cases:
(i) $v_k>v_n+1$,
(ii) $v_k=v_n+1$,
(iii) $v_k=v_n$,
(iv) $v_k<v_n$.

\vspace{.1in}
\noindent
(i) $v_k>v_n+1$. Let $m$ be an integer such that
$0\leq m< v_k-(v_n+1)$.
It implies $0< m+1< v_k-v_n$.
Due to the assumption $P(\gamma)$,
we have the inequality
$\deg_t(\xi(\vecv)^{\tau_{n,k}(\vecv)+(m+1)\vece_k})\leq
 \co(\vecv)-2$.
If $v_j<v_n<v_k$, then we have $j\neq k$.
Since
$q_j(\tau_{n,k}(\vecv)+(m+1)\vece_k)=v_j<v_n=q_j(\tau_{n,j}(\vecv))$,
we have
$\tau_{n,k}(\vecv)+(m+1)\vece_k\not\geqzero \tau_{n,j}(\vecv)$.
Thus we have the inequality
$\deg_t\xi(\tau_{n,j}\vecv)^{\tau_{n,k}\vecv+(m+1)\cdot \vece_k}
 <\co(\vecv)$.
Due to (\ref{eq;1.9.30}),
we have the following inequality:
\[
 \deg_t\xi(\vecv+\vece_n)^{\tau_{n\,k}(\vecv+\vece_n)+m\cdot\vece_k}
\leq
 \co(\vecv)-1\leq \co(\vecv+\vece_n)-2.
\]
Thus we are done in this case.

\vspace{.1in}
\noindent
(ii) $v_k=v_n+1$.
In this case, the equality (\ref{eq;1.9.31}) is a consequence of 
Lemma \ref{lem;1.9.20}.

\vspace{.1in}
\noindent
(iii) $v_k=v_n$.
When $v_j<v_n$, we have $j\neq k$.
Thus
we have $\tau_{n,j}(\vecv)\not\leqzero \vecv+(m+1)\vece_k$,
since $q_j(\tau_{n,j}\vecv)=v_n>v_j=q_j(\vecv+(m+1)\cdot\vece_k)$.
Hence we have the inequality
$\deg_t\xi(\tau_{n,j}\vecv)^{\vecv+(m+1)\cdot\vece_k}
\leq \co(\vecv)-1$.
On the other hand, we have the equalities
$\deg_t\xi(\vecv)^{\vecv+(m+1)\vece_k}=\co(\vecv)$
and
$\omega_n^{\vecv+(m+1)\vece_k}-\omega_n^{\vecv}=\Delta_{n,k}$.
Thus we obtain the following:
\begin{multline}
 \big[\xi(\vecv+\vece_n)^{\tau_{n,k}(\vecv+\vece_n)+m\cdot\vece_k}\big]
=[\omega_n^{\vecv+(m+1)\cdot\vece_k}-\omega_n^{\vecv}]
\cdot
 [\xi(\vecv)^{\vecv+(m+1)\vece_k}]\\
=\Delta_{n,k}\cdot
 \prod_{j\neq k}\prod_{i=0}^{v_{j}-1}(t_{v_j}-t_i)
\times
 \prod_{i=0}^{v_k-1}(t_{v_k+1+m}-t_i)
=\Delta_{n,k}\cdot
 \prod_{j=1}^{n-1}\prod_{i=0}^{v_{j}-1}(t_{v_j}-t_i)
\times
 \prod_{i=0}^{v_n-1}(t_{v_n+1+m}-t_i).
\end{multline}

\noindent
(iv) $v_k<v_n$.
We have the equality
$[\omega_n^{\tau_{n,k}\vecv+(m+1)\cdot \vece_k}-
  \omega_n^{\vecv}]=(t_{v_k}-t_{v_n})$.
We also have the following:
\[
 [\xi(\vecv)^{\tau_{n,k}\vecv+(m+1)\cdot\vece_k}]
=\Delta_{k,n}\cdot
 \prod_{j=1}^{n-1}\prod_{i=0}^{v_j-1}(t_{v_j}-t_i)
\times
\prod_{\substack{0\leq i\leq v_n-1\\ i\neq v_k}}
 (t_{v_n+m+1}-t_i).
\]
In particular,
the total degree of the product is $\co(\vecv)$.

If $k\neq j$,
we have the inequality
$\deg_t\xi(\tau_{k\,n}\vecv)^{\tau_{n,k}\vecv+(m+1)\cdot\vece_k}
<\co(\vecv)$.
If $k=j$, we have the following equality:
\[
 [\xi(\tau_{n,k}\vecv)^{\tau_{n,k}\vecv+(m+1)\vece_k}]=
 \prod_{j=1}^{n-1}\prod_{i=0}^{v_j-1}(t_{v_j}-t_i)
\times
 \prod_{i=0}^{v_n-1}(t_{v_n+m+1}-t_i).
\]
Thus we obtain the following by a direct calculation:
\[
  [\xi(\vecv+\vece_k)^{\tau_{n,k}(\vecv+\vece_k)+m\cdot\vece_k}]=
\Delta_{k\,n}\cdot
 \prod_{j=1}^{n-1}\prod_{i=0}^{v_j-1}
(t_{v_j}-t_i)
\times
 \prod_{\substack{0\leq i\leq v_n\\ i\neq v_k}}
(t_{v_n+1+m}-t_i).
\]
Hence the proof of Lemma \ref{lem;1.9.35} is completed.
\hfill\qed

\vspace{.in}
Thus the induction of the proof can proceed,
namely, the proof of Theorem \ref{thm;1.9.40} is completed.
\hfill\qed

\begin{cor}
We have the following equality:
$\xi(\vece_i)=\omega_i-t_0-p_i^{\ast}c_1(L_0)
  +\sum_{k<i}\Delta_{k,i}$.
Thus we have the following equalities for $\vecv\in [0,r-1]^n$:
\[
 \xi(\vece_i)^{\vecv}=
 (t_{v_i}-t_0)+\sum _{\substack{k<i\\v_i\neq v_k}}\Delta_{i\,k}
+p_i^{\ast}\bigl(c_1(L_{v_i})-c_1(L_0)\bigr).
\]
\hfill\qed
\end{cor}

\subsection{Some corollaries}

\subsubsection{The product in the associated graded ring}

Recall that we have the filtration
$\{F_lH^{\ast}(\Filt(E,n))\,|\,l\in\seisuuplus\}$.
Let denote the associated graded space
by $Gr H^{\ast}(\Filt(E,n))$.
\begin{cor} \label{cor;1.9.50}
For any $\vecv\in[0,r-1]^n$,
$\omega_n\cdot \xi(\vecv)$ is contained in
$F_{\co(\vecv)+1}H^{\ast}(\Filt(E,n))$.
Moreover we have the following equality in
 $Gr_{\co(\vecv)+1}H^{\ast}(\Filt(E,n))$:
\[
 \omega_n\cdot \xi(\vecv)\equiv \xi(\vecv+\vece_n).
\]
Here we use the notation $\equiv$ to denote the equality
in $Gr_{\co(\vecv)+1}H^{\ast}(\Filt(E,n))$.
\end{cor}
\pf
Clear from Theorem \ref{thm;1.9.40}.
\hfill\qed

\begin{cor} \label{cor;1.13.40}
 Let $\vecv$ be an element of $[0,r-1]^n$.
 Then we have the following equalities in
the associated graded ring
$Gr_{\co(\vecv)} H^{\ast}(\Filt(E,n))$:
 \[
  \xi(\vecv)\equiv \prod_{i=1}^n \xi(\vece_i)^{v_i}
 \equiv \prod_{i=1}^n\omega_i^{v_i}.
 \]
Let $\vecv$ and $\vecv'$ be elements of $[0,r-1]^n$.
Then we have 
$\xi(\vecv)\cdot\xi(\vecv')\equiv\xi(\vecv+\vecv')$
in $Gr_{\co(\vecv+\vecv')}H^{\ast}(\Filt(E,n))$.
\end{cor}
\pf
We only have to use Corollary \ref{cor;1.9.50} inductively.
\hfill\qed

\begin{cor}
The ring $Gr H^{\ast}(\Filt(E,n))$ is isomorphic to
the following:
\[
 H^{\ast}(C^n)[\omega_1,\ldots,\omega_n]/(\omega_i^r=0).
\]
If we are given $L_i$ $(i=0,1,\ldots,)$ and if we put
$E_{\infty}=\bigoplus_i E_i$,
then the ring $Gr H^{\ast}(\Filt(E_{\infty},n))$
is isomorphic to the polynomial ring
over $H^{\ast}(C^n)$ with $n$ variables.
\hfill\qed
\end{cor}

\subsubsection{Some other equalities}

Let $\vecv$ be an element of $[0,r-1]^n$ such that $v_n\neq 0$.
Let $m$ be an integer such that $0\leq m<n$.
We put as follows:
\[
 H_{m,n}(\vecv)=
 \left\{
 \begin{array}{ll}
 0 & (v_m\geq v_n)\\
 1 & (v_m<v_n).
 \end{array}
 \right.
\]
\begin{thm} \label{thm;1.9.60}
In the equivariant cohomology ring, we have the following equality:
\begin{equation}
 \Bigl(
 \omega_m-\sum_{\substack{k<m,\\ v_k=v_m}}
 \Delta_{k,m}-p_m^{\ast}L_{v_m}-t_{v_m}
 \Bigr)\cdot \xi(\vecv)
=\xi(\vecv+\vece_m)
+H_{m,n}(\vecv)\cdot \Delta_{m,n}\xi(\tau_{m,n}\vecv).
\end{equation}
As a corollary, we obtain the following equality in the non-equivariant
cohomology ring:
\begin{equation} \label{eq;1.9.70}
 \Bigl(
 \omega_m-\sum_{\substack{k<m,\\v_k=v_m}}\Delta_{k,m}-p_m^{\ast}L_{v_m}
 \Bigr)\cdot \xi(\vecv)
=\xi(\vecv+\vece_m)
+H_{m,n}(\vecv)\cdot \Delta_{m,n}\xi(\tau_{m,n}\vecv).
\end{equation}
\end{thm}
\pf
We put $X=(\omega_m-\omega_m^{\vecv})\cdot \xi(\vecv)$.
Due to an argument similar to the proof of Lemma \ref{lem;1.9.10},
we obtain the following description of $X$:
\[
 X=\xi(\vecv+\vece_m)+
\sum_{v_k\neq v_m}A_k\cdot \xi(\tau_{k\,m}\vecv).
\]
Here $A_k$ are elements of $H^2(C^n)$.

\begin{lem}
Assume $\deg_t\xi(\vecv)^{\tau_{m,k}(\vecv)}<\co(\vecv)-1$.
Then we have $A_k=0$.
\end{lem}
\pf
By the assumption, we have $\deg_t X^{\tau_{m,k}}(\vecv)<\co(\vecv)$.
We have $\xi(\vecv+\vece_m)^{\tau_{m,k}(\vecv)}=0$.
We have
$\deg_t \xi(\tau_{j\,m}\vecv)^{\tau_{k\,m}(\vecv)}<\co(\vecv)$
for any $j\neq k$, when $v_j\neq v_m$ and $v_k\neq v_m$.
On the other hand, we have the following:
\[
 [\xi(\tau_{m,k}\vecv)^{\tau_{m,k}\vecv}]
=\prod_{j=1}^n\prod_{i=0}^{v_j-1}
 (t_{v_j}-t_i).
\]
In particular, we have
$\deg_t\xi(\tau_{m,k}\vecv)^{\tau_{m,k}\vecv}
=\co(\vecv)$.
Thus we obtain $A_k=0$.
\hfill\qed

\vspace{.1in}

Assume that $k=n$.
If $v_m>v_n$, then we have
$\deg_t\xi(\vecv)^{\tau_{m,n}}<\co(\vecv)-1$.
Hence we obtain $A_n=0$ in this case.
If $v_m<v_n$, then we have the following:
\[
 [\xi(\vecv)^{\tau_{m,n}(\vecv)}]
=\Delta_{m,n}\cdot
 \prod_{j=1}^{n-1}\prod_{i=0}^{v_j-1}
 (t_{v_j}-t_i)
\times
 \prod_{\substack{0\leq i\leq v_n-1\\i\neq v_m}}(t_{v_n}-t_i).
\]
We also have
$[\omega_m^{\tau_{m,n}(\vecv)}-\omega_m^{\vecv}]
=t_{v_n}-t_{v_m}$.
Thus the total degree of the product is $\co(\vecv)$.
By comparing the top terms,
we obtain $A_n=\Delta_{m,n}$ in this case.

Then Theorem \ref{thm;1.9.60} is a corollary of the next lemma.
\begin{lem}
Let consider the case $k\neq n$.
Assume $v_k\neq v_n$.
We have the inequality
$\deg_t\xi(\vecv)^{\tau_{m,k}\vecv}<\co(\vecv)-1$.
As a corollary, we obtain $A_k=0$ if $k\neq n$.
\end{lem}
\pf
We already have the following equality:
\[
 \xi(\vecv)=
 (\omega_n-\omega_n^{\vecv-\vece_n})\cdot 
 \xi(\vecv-\vece_n)
+\sum_{\substack{j<n,\\ v_j<v_{n}-1}}
 \Delta_{j,n}\cdot \xi(\tau_{j,n}(\vecv-\vece_n)).
\]
We have the following:
\[
 \deg_t(\omega_n^{\tau_{k,m}(\vecv-\vece_n)}-\omega_n^{\vecv-\vece_n}) 
=0,
\quad
 \deg_t\xi(\vecv-\vece_n)^{\tau_{k,m}(\vecv-\vece_n)}
 \leq \co(\vecv-\vece_n)-1
 <\co(\vecv)-1.
\]
We also have the following inequality:
\[
 \deg_t\xi(\tau_{j,n}(\vecv-\vece_n))^{\tau_{k,m}\vecv}
\leq
 \co(\vecv-\vece_n)-1<\co(\vecv)-1.
\]
In all, we obtain the inequality desired.
\hfill\qed

\begin{rem}
Clearly the formulas $(\ref{eq;1.9.71})$ and $(\ref{eq;1.9.70})$
provide the method to describe $\xi(\vecv)$
in terms of $\omega_i$ $(i=1,\ldots,n)$ explicitly.
\hfill\qed
\end{rem}

%% file: 12.4.tex
\section{The structure of the cohomology ring of filt schemes}

\subsection{Preliminary}

Let $n$ be a positive integer.
We have the left action of $\gbigs_n$ on $\seisuuplus$
by the permutation of the components:
\[
 \sigma(a_1,\ldots,a_n):=
 (a_{\sigma^{-1}(1)},\ldots,a_{\sigma^{-1}(n)}).
\]
Let $\vece_i$ $(i=1,\ldots,n)$ be the canonical base of $\seisuu^n$.
Then $\sigma(\vece_i)=\vece_{\sigma(i)}$.
For an element $\vecv\in\seisuu^n$,
we put $St(\vecv):=\{\sigma\in\gbigs_n\,|\,\sigma\vecv=\vecv\}$.

Similarly we have the left action of $\gbigs_n$
on $C^n$ by the permutation of the components:
$\sigma(x_1,\ldots,x_n)=(x_{\sigma^{-1}(1)},\ldots,x_{\sigma^{-1}(n)})$.
Thus we have the left action of $\gbigs_n$
on $C^n\times \seisuuplus^n$.
We have the natural identification $F(\vecv)=C^n$.
Thus we have the $\gbigs_n$-action on
$\coprod_{\vecv\in [0,r-1]^n} F(\vecv)$.

For an element $\sigma\in\gbigs_n$,
we obtain the automorphism
$(\sigma^{-1})^{\ast}$
of $\bigoplus_{\vecv\in [0,r-1]}H^{\ast}(F(\vecv))$.
We denote $(\sigma^{-1})^{\ast}(a)=\sigma\cdot a$.
Thus we obtain the $\gbigs_n$-action on
$\bigoplus_{\vecv\in [0,r-1]}H^{\ast}(F(\vecv))$.

Let $E=\bigoplus_{\alpha=0}^{r-1}L_{\alpha}$
be a direct sum of line bundles on $C$.
We introduce the $\gbigs_n$-action on $H^{\ast}(\Filt(E,n))$
as follows:
We have the decomposition
$H^{\ast}(\Filt(E,n))=
 \bigoplus_{\vecv\in[0,r-1]^n} H^{\ast}(F(\vecv))\cdot\xi(\vecv)$.
Then we put as follows:
\[
 \sigma\Bigl(\xi(\vecv)\cdot a
 \Bigr)=
 \xi(\sigma(\vecv))\cdot \sigma(a).
\]
It is clear that the restriction
$H^{\ast}(\Filt(E,n))\lrarr
 \bigoplus_{\vecv\in[0,r-1]^n}H^{\ast}(F(\vecv))$
is equivariant with respect to the $\gbigs_n$-actions above.

We can also consider similar $\gbigs_n$-actions
on the equivariant cohomology groups
$H^{\ast}_{G_m^r}(\Filt(E,n))$
and
$\bigoplus_{\vecv\in[0,r-1]^n}
  H^{\ast}_{G_m^r}(F(\vecv))$.

\begin{lem}
The $\gbigs_n$-action preserves the filtration
$F_hH^{\ast}(\Filt(E,n))$.
Thus we obtain the $\gbigs_n$-action on 
the associated graded vector space $\Gr H^{\ast}(\Filt(E,n))$.
\end{lem}
\pf
Clear from our constructions.
\hfill\qed

The relation of the product and the $\gbigs_n$-action
on $H^{\ast}(\Filt(E,n))$ is not clear for the author,
at the moment.
However the following lemma is easy to see.
\begin{lem}
The $\gbigs_n$-action preserves the product of
$\Gr H^{\ast}(\Filt(E,n))$.
\end{lem}
\pf
Use Corollary \ref{cor;1.13.40}.
\hfill\qed

\subsection{Rough structure}
\label{subsection;1.13.1}

\subsubsection{The structure of the associated graded ring}
\label{subsubsection;1.13.5}

Let $\vecv$ be an element of $[0,r-1]^n$.
We have the normalization $\Nor(\vecv)\in \nbigb(n,r)$.
Then we have the morphism $F(\vecv)\lrarr F^Q(\Nor(\vecv))$.
Thus we obtain the following morphism
for $\vecu\in \nbigb(n,r)$:
\[
 \coprod_{\Nor(\vecv)=\vecu}
 F(\vecv)
\lrarr
 F^Q(\vecu).
\]
We already have the $\gbigs_n$-action on the left hand side.
On the other hand, we have the trivial action
on the right hand side.
The morphism is equivariant.

In particular, we have the morphism
$H^{\ast}(F^Q(\vecu))\lrarr H^{\ast}(F(\vecu))$.
Here we regard $\vecu\in\nbigb(n,r)$ as the element of $[0,r-1]^n$
by the natural inclusion.
Then $H^{\ast}(F^Q(\vecu))$
is isomorphic to the $St(\vecu)$-invariant part
of $H^{\ast}(F(\vecu))$.
We identify them by the morphism.

If $\Nor(\vecv)=\vecu$,
then we have an element $\sigma(\vecu)=\vecv$.
The morphism
$H^{\ast}(F^Q(\vecu))\lrarr H^{\ast}(F(\vecv))$ is 
identified with $a\longmapsto \sigma\cdot a$.
It is same as the following:
\[
 a\longmapsto \frac{1}{|St(\vecu)|} 
 \sum_{\sigma\in\gbigs_n,\sigma(\vecu)=\vecv}
 \sigma\cdot a.
\]

\begin{lem} \label{lem;1.13.50}
We have the following equality:
\[
 \Psi^{\ast}\xi^Q(\vecu;a)=
 \frac{1}{|St(\vecu)|}
\sum_{\sigma\in\gbigs_n}
 \sigma\Bigl(
 \xi(\vecu)\cdot a
 \Bigr).
\]
\end{lem}
\pf
We use a standard argument to use the equivariant cohomology ring.
We only indicate the outline.
By a geometric consideration,
we can show that $\deg_t\xi^Q(\vecu;a)\leq \co(\vecu)$.
We describe $\Psi^{\ast}\xi^Q(\vecu;a)$
as $\sum_{\vecv\in [0,r-1]^n} A(\vecv)\cdot\xi(\vecv)$.
\begin{lem}
$A(\vecv)=0$ unless $\Nor(\vecv)\geq_0\vecu$.
\end{lem}
\pf
Similar to Lemma \ref{lem;1.9.10}.
\hfill\qed

\begin{lem}
Let $\vecv$ be a permutation of $\vecu$.
Then we have
 $A(\vecv)=|St(\vecu)|^{-1}\sum_{\sigma\vecu=\vecv}\sigma(a)$.
\end{lem}
\pf
We only have to consider the restriction
of $\Psi^{\ast}\xi^Q(\vecu;a)$
to $F(\vecv)$.
\hfill\qed

\vspace{.1in}
We put
$Y:=
 \Psi^{\ast}\xi^Q(\vecu;a)
 -|St(\vecu)|^{-1}
 \sum_{\sigma\in\gbigs_n}\sigma\bigl(\xi(\vecu)\cdot a\bigr)$.
We have $\deg_tY^{\vecw}\leqzero\co(\vecu)$ for any $\vecw$.
We also have $Y^{\vecw}=0$ unless $N(\vecw)>_0\vecu$.
Then we can conclude that $Y=0$.
Hence the proof of Lemma \ref{lem;1.13.50} is completed.
\hfill\qed

We have some easy corollaries.
\begin{cor}
\mbox{{}}
\begin{itemize}
\item
The image $\Psi^{\ast}(H^{\ast}(\Quot(E,n)))$
is same as the $\gbigs_n$-invariant part
of $H^{\ast}(\Filt(E,n))$.
\item
The morphism $\Psi^{\ast}$ preserves the filtrations,
that is,
$\Psi^{\ast}(F_hH^{\ast}(\Quot(E,n)))$
is contained in $F_hH^{\ast}(\Filt(E,n))$.
Moreover we have the following:
\[
\Psi^{\ast}\bigl(F_hH^{\ast}(\Quot(E,n))\bigr)
 \cap F_{h-1}H^{\ast}(\Filt(E,n))
=\Psi^{\ast}\bigl(F_{h-1}H^{\ast}(\Quot(E,n))\bigr).
\]
\item
The filtration $F_hH^{\ast}(\Quot(E,n))$ is compatible
with the product of $H^{\ast}(\Quot(E,n))$,
that is,
$F_h\cdot F_k\subset F_{h+k}$.

\item
Thus we obtain the ring morphism
$\Gr(\Psi^{\ast}):\Gr H^{\ast}(\Quot(E,n))
 \lrarr \Gr H^{\ast}(\Filt(E,n))$.
The image of $\Gr(\Psi^{\ast})$
is same as the $\gbigs_n$-invariant part of
$\Gr H^{\ast}(\Filt(E,n))$.
\hfill\qed
\end{itemize}
\end{cor}

\begin{cor}
In $\Gr H^{\ast}(\Quot(E,n))$,
we have the following equality:
\[
 \Psi^{\ast}\xi^Q(\vecu;a)
\equiv
 \frac{1}{|St(\vecu)|}
 \sum_{\sigma\in\gbigs_n}
 \sigma\Bigl(
 \prod_{i=1}^n\omega_i^{u_i}\cdot a
 \Bigr).
\]
\hfill\qed
\end{cor}

\subsubsection{Another filtration}

We have the lexicographic order on $\seisuuplus^{n+1}$.
We have the natural injection
$\seisuuplus^{n}\lrarr \seisuuplus^{n+1}$
given by $\vecv\longmapsto (\co(\vecv),\vecv)$.
Then we obtain the induced total order on $\seisuuplus^n$.
We denote the induced order by $\leqiti$.

We put as follows:
\[
 F^{(1)}_{\vecv}H^{\ast}(\Filt(E,n))
=\bigoplus_{\vecu\leqiti\vecv}
 H^{\ast}(F(\vecu))\cdot\xi(\vecu).
\]
We use the notation $\equiv^{(1)}$ to denote the equality
in the associated graded vector space
$\Gr^{(1)}H^{\ast}(\Filt(E,n))$.
The product
$\xi(\vecv)\cdot \xi(\vece_i)$
is contained in $F_{\vecv+\vece_i}^{(1)}$,
and same as $\xi(\vecv+\vece_i)$
in $Gr^{(1)}_{\vecv+\vece_i}$.
Thus we can show that
$F^{(1)}_{\vecv}\cdot F^{(1)}_{\vecu}\subset F^{(1)}_{\vecv+\vecu}$,
and $\xi(\vecv)\cdot \xi(\vecu)\equiv^{(1)}\xi(\vecv+\vecu)$
in $\Gr^{(1)}_{\vecv+\vecu}$.

Let $\vecv$ be an element of $\nbigb(n,r)$.
Then $\sigma(\vecv)<_1\vecv$ if $\sigma\neq 1$.
Hence $\xi(\sigma\vecv)$ is $0$ in $\Gr^{(1)}_{\vecv}$.
In particular,
we obtain the following equality in $\Gr^{(1)}_{\vecv}$
for $\vecv\in\nbigb(n,r)$:
\[
 \frac{1}{|St(\vecv)|}
\sum_{\sigma\in\gbigs_n}
 \sigma(\xi(\vecv)\cdot a)
\equiv^{(1)}
 \xi(\vecv)\cdot a.
\]

\subsubsection{Some generators}

Recall that $\vecc_l$ denotes $\sum_{j\leq l}\vece_j$.
Then $F^Q(\vecc_l)$ is isomorphic to $C^{(l)}\times C^{(n-l)}$.
Thus we have the natural inclusion
$f:H^{\ast}(C^{(l)})\lrarr H^{\ast}(F^Q(\vecc_l))$,
induced by the projection.
We also have the morphism
$g:H^{\ast}(C^{(n)})\lrarr H^{\ast}(C^{(l)})\otimes
H^{\ast}(C^{(n-l)})$
induced by the projection $C^{(n)}\lrarr C^{(l)}\times C^{(n-l)}$.
Then we obtain the morphism
$\phi:
 H^{\ast}(C^{(l)})\otimes
 H^{\ast}(C^{(n)})
\lrarr
 H^{\ast}(C^{(l)})\otimes H^{\ast}(C^{(n-l)})$,
given by $\phi(a\otimes b)=f(a)\cdot g(b)$.

\begin{lem}
The morphism $\phi$ is surjective.
\end{lem}
\pf
We have the filtration on
the $H^{\ast}(C^{(k)})$ and the $H^{\ast}(C^{(k-j)})$
defined as follows:
We denote the subspace generated by the elements
of the form
$\sym(
 \gamma_1\otimes\cdots\otimes\gamma_a
 \otimes 1\otimes\cdots\otimes 1),\,
 (\gamma_i\in H^{\ast}(C))$
by $G_a$.
Here $\sym$ denotes the symmetrization operator.
Then the family $\{G_a\}$ gives the filtration.
Moreover the family
$\big\{H^{\ast}(C^{(j)})\otimes G_a(H^{\ast}(C^{(k-j)}))
 \big\}$
gives the filtration on the space
$H^{\ast}(C^{(j)})\otimes H^{\ast}(C^{(k-j)})$.

\vspace{.1in}

\noindent 
{\bf Claim}\quad\quad
The space $H^{\ast}(C^{(j)})\otimes G_aH^{\ast}(C^{k-j})$
is contained in the image $\Image(\phi)$.

\vspace{.1in}

We prove the claim by an induction on the number $a$.
Assume that the claims for any $a'<a$ hold.
For any
$\gamma=\gamma_1\otimes\cdots\otimes\gamma_a\otimes 1^{\otimes\,k-j-a}
\in H^{\ast}(C)^{\otimes k-j}$,
we have the description of
$\sym_{k}(1^{\otimes\,j}\otimes\gamma)$
as follows:
\[
 \sym_k(1^{\otimes\,j}\otimes\gamma)=
 1\otimes \sym_{k-j}(\gamma)+
 \sum_{a'<a} \sum _i\ \sym_{j}(c_{a',i})\otimes \sym_{k-j}(c'_{a',i}),
\quad
 \sym_{k-j}(c'_{a',i})\in G_{a'}H^{\ast}(C^{(k-j)}).
\]
Due to the equality,
the following holds
for any $\gamma'\in H^{\ast}(C^{(j)})$:
\[
\gamma' \otimes \sym_{k-j}(\gamma)
=\gamma'\otimes 1\cdot
 \Bigl(\sym_{k}(1^{\otimes\,j}\otimes\gamma)\Bigr)
-\sum_{a'<a}\sum_{i}(\sym_{j}(c_{a',i})\gamma')\otimes
\ sym_{k-j}(c'_{a',i})
\]
which belongs to the subgroup
$Im(\phi)+H^{\ast}(C^{(j)})\otimes G_{a-1}H^{\ast}(C^{(k-j)})$.
Thus we are done.
\hfill\qed

For any $\vech\in\seisuuplus^h$,
we have the following morphism induced by
$\phi$ and the identify of $C^{(\vech)}$:
\[
 H^{\ast}(C^{(j)})\otimes H^{\ast}(C^{(k)})
 \otimes H^{\ast}(C^{(\vech)})
\lrarr
 H^{\ast}(C^{(j)})\otimes H^{\ast}(C^{(k-j)})
 \otimes H^{\ast}(C^{(\vech)}).
\]
We denote it also by $\phi$.
\begin{lem} \label{lem;1.13.60}
The morphism $\phi$ is surjective.
\hfill\qed
\end{lem}

Let $L_{\alpha}$ $(\alpha=0,1,2,\ldots)$ be line bundles over $C$.
We put $E_{\infty}=\bigoplus_{\alpha\geq 0}L_{\alpha}$.
We obtain the infinite quot scheme $\Quot(E_{\infty},n)$.
Since $F(\vecc_l)$ is isomorphic to $C^{(l)}\times C^{(n-l)}$,
we have the natural inclusion
$H^{\ast}(C^{(l)})\subset H^{\ast}(F(\vecc_l))$.
\begin{prop} \label{prop;1.13.10}
The following elements generate $H^{\ast}(\Quot(E_{\infty},n))$
over the ring $H^{\ast}(C^{(n)})$:
\[
 S:=
 \big\{\xi^Q(\vecc_l;a)\,|\,l=1,\ldots,n,\,\,
 a\in H^{\ast}(C^{(l)})\subset H^{\ast}(F(\vecc_l))\big\}.
\]
\end{prop}
\pf
Let $\nbiga$ denote the subalgebra generated by $S$ over $H^{\ast}(C^{(n)})$.
We have to show that $\nbiga$ is same as $H^{\ast}(\Quot(E_{\infty},n))$.
Note the following equality in $\Gr_{\vecu}^{(1)}$:
\[
 \Psi^{\ast}\xi^Q(\vecu;a)
\equiv^{(1)}
 \xi(\vecu;a).
\]
We use the following obvious lemma.
\begin{lem}
Let $\vecu$ be an element of $\nbigb(n,r)$.
Let $k$ be a positive integer such that 
$k<\min\{i\,|\,u_i>u_{i+1}\}$.
Then we have the following equality in $\Gr^{(1)}_{\vecu+\vecc_k}$:
\[
 \Psi^{\ast}\xi^Q(\vecu;a)\cdot
 \Psi^{\ast}\xi(\vecc_l;b)
\equiv^{(1)}
 \Psi^{\ast}\xi^Q(\vecu+\vecc_l;\phi(a\otimes b))
\]
\hfill\qed
\end{lem}

For any element $f$ of $H^{\ast}(\Quot(E_{\infty},n))$,
it has a description as follows:
\[
 f= \sum_{\vecu\in\nbigb(n,r)}\sigma\xi^Q(\vecu;a(\vecu)).
\]
We put $d(f):=\max\{\vecu\,|\,a(\vecu)\neq 0\}$.
Here the maximum is taken with respect to
the total order $\leqiti$ on $\nbigb(n,r)$.

Assume that $H^{\ast}(\Quot(E_{\infty},n))-\nbiga\neq \emptyset$,
and we will derive a contradiction.
Let $f_0$ be an element of $H^{\ast}(\Quot(E_{\infty},n))-\nbiga$
satisfying
$d(f_0)=\min\bigl\{d(f)\,|\,f\in H^{\ast}(\Quot(E_{\infty},n))-\nbiga\bigr\}$.
We put $\vecu_0:=d(f_0)$.

We put $k:=\min\{i\,|\,u_{0,i}>u_{0,i+1}\}$.
Then we have the decomposition $\vecu_0=\vecu_1+\vecc_k$.
Note that $\vecu_1<\vecu_0$.
Thus we only have to show that
there exist $a\in H^{\ast}(F(\vecu_1))$
and $b\in H^{\ast}(C^{(k)})$
satisfying the following:
\[
 d(f_0-\xi^Q(\vecu_1;a)\cdot \xi^Q(\vecc_k;b))
<d(f_0).
\]
In $Gr^{(1)}_{\vecu_0}$,
we have the following:
\[
 f_0\equiv^{(1)}\xi^Q(\vecu_0;a(\vecu_0)),
 \quad
 \xi^Q(\vecu_1;a)\cdot \xi^Q(\vecc_k;b)
\equiv_1\xi^Q(\vecu+\vecc_k;\phi(a\otimes b)).
\]
Thus we obtain the result,
by using Lemma \ref{lem;1.13.60}.
\hfill\qed

\vspace{.1in}

Let $k$ be a positive integer.
Then we have the push forward
$q_{\ast}:H^{\ast}(C^{(k-1)})\otimes H^{\ast}(C)\lrarr H^{\ast}(C^{(k)})$.
Let $\phi_2$ denote the composite of the following morphisms:
\[
\begin{CD}
 H^{\ast}(C^{(k-1)})\otimes H^{\ast}(C^{(l)})\otimes H^{\ast}(C)
 @>{id\otimes \phi}>>
 H^{\ast}(C^{(k-1)})\otimes H^{\ast}(C)\otimes H^{\ast}(C^{(l-1)})
 @>{q_{\ast}\otimes id}>>
 H^{\ast}(C^{(k)})\otimes H^{\ast}(C^{(l-1)}).
\end{CD}
\]
Clearly $\phi_2$ is surjective.
In general we have the morphism induced by $\phi_2$:
\[
 H^{\ast}(C^{(\vech)})\otimes
 H^{\ast}(C^{(k-1)})\otimes H^{\ast}(C^{(l)})
\otimes H^{\ast}(C)
\lrarr H^{\ast}(C^{(\vech)})
\otimes H^{\ast}(C^{(k)})\otimes H^{\ast}(C^{(l-1)}).
\]
The morphism is also denoted by $\phi_2$.
It is clearly surjective.

\begin{prop} \label{prop;1.13.11}
The following set generates the ring $H^{\ast}(\Quot(\infty,n))$
over $H^{\ast}(C^{(n)})$:
\[
 S_2=
 \{\xi^Q(l\cdot\vece_1,p^{\ast}a)\,|\,
 a\in H^{\ast}(C),\,\,l=1,\ldots,n\}.
\]
\end{prop}
\pf
Let $\nbiga_2$ denote the subring generated by $S_2$
over $H^{\ast}(C^{(n)})$.
We only have to prove that $\nbiga_2$
contains $S$ in Proposition \ref{prop;1.13.10}.
For $k\leq n$, we have the natural inclusion
$\nbigb(k,r)\lrarr\nbigb(n,r)$.
We regard $\nbigb(k,r)$ as the subset of $\nbigb(n,r)$.

We use an induction.
The following claim is called $P(l,k)$:
\begin{quote}
$P(l,k)$:
For any element $\vecv\in\nbigb(k,r)$ such that
$\co(\vecv)\leq l$,
the element $\xi^Q(\vecv;a)$ is contained in $\nbiga_2$
for any $a\in H^{\ast}(F(\vecv))$.
\end{quote}
We only have to prove $P(n,n)$.
The claims $P(l,1)$ $(1\leq l\leq n)$ are clear.
Thus the proposition is reduced to the following lemma.
\begin{lem}
$P(l-1,l-1)+P(l,k-1)\Longrightarrow P(l,k)$
for $k\leq l$.
\end{lem}
\pf
Let $\vecu$ be an element of $\nbigb(k,r)-\nbigb(k-1,r)$
such that $\co(\vecu)\leq l$.
The $k$-th component is denoted by $u_k$.
We put $\vecu_0:=\vecu-u_k\cdot \vece_k$.
Then we have $\co(\vecu_0)\leq l-u_k<l$.
We have the following equality in
$Gr_{\co(\vecu)} H^{\ast}(\Quot(\infty,n))$:
\[
 \xi^Q(\vecu_0;a)\cdot \xi^Q(u_k\cdot\vece_k;b)
\equiv
 \xi^Q(\vecu;\phi_2(a\otimes b))
+\sum_{\vecw\in\nbigb(k-1,r)}\xi^Q(\vecw;c(\vecw)),
\quad
 (c(\vecw)\in H^{\ast}(F(\vecw))).
\]
Thus the lemma follows from the surjectivity of $\phi_2$.
\hfill\qed

Hence the proof of the proposition \ref{prop;1.13.11}
is completed.
\hfill\qed

\vspace{.1in}

Let $E$ be a locally free sheaf,
which is not necessarily isomorphic to a direct sum
of line bundles.
We have the universal quotient $\nbigfu$
defined over $\Quot(E,n)\times C$.
We have the $i$-th Chern class $c_i(\nbigfu)$.
It induces the correspondence map
$\phi_i:H^{\ast}(C)[-2(i-1)]\lrarr H^{\ast}(\Quot(E,n))$.

\begin{thm} \label{thm;1.13.12}
The following set generates $H^{\ast}(\Quot(E,n))$:
\[
 \bigcup_{i=1}^{n+1}\Image(\phi_i).
\]
\end{thm}
\pf
Because any vector bundle $E$
on a smooth projective curve can be deformed into such a bundle,
and because the topology of the quot scheme is invariant
under deformations of vector bundles,
we only have to prove the claim for the vector bundle
of the form $\bigoplus_{\alpha=0}^{r-1} L_{\alpha}$.
We have the morphism
$\Psi\times id:\Filt(E,n)\times C\lrarr \Quot(E,n)\times C$.
We have the filtrations
$\nbigfu_n\lrarr\nbigfu_{n-1}\lrarr\cdots\lrarr\nbigfu_1$
over $\Filt(E,n)\times C$,
and $(\Psi\times id)^{\ast}\nbigfu=\nbigfu_n$.
We put $\nbiggu_j:=\ker(\nbigfu_j\lrarr\nbigfu_{j-1})$.
We have already known the following:
\[
 c_{\ast}(\nbiggu_j)=
 1+\Delta_{j\,n+1}\sum _{h=0}^{\infty}(-\omega_i)^h
+\Delta_{j,n+1}^2(1+\omega_j)^{-2}.
\]
Hence we obtain the following:
\begin{multline}
 (\Psi\times id)^{\ast}
 c_{\ast}(\nbigf)
=\prod_{i}
 \Bigl(
 1+\Delta_{j,n+1}\cdot\sum_{h=0}^{\infty}(-\omega_j)^h
 +\Delta_{j,n+1}^2\cdot(1+\omega_j)^{-2}
 \Bigr)\\
=
1+\sum_{h=0}^{\infty}
  \sum_{j=1}^n\Delta_{j,n+1}(-\omega_j)^h
+A.
\end{multline}
Here $A$ is contained in the ideal generated by
$H^4(C^n)$.

We have $\phi_1(a)=\sum p_i^{\ast}(a)\in H^{\ast}(C^{(n)})$.
For $i>1$, we have the equality
$\phi_i(a)\equiv \xi((i-1)\vece_1;a)$
in $Gr_{i-1}H^{\ast}(\Filt(E,n))$.
Hence Theorem \ref{thm;1.13.12} follows from
Proposition \ref{prop;1.13.11}.
\hfill\qed

%% file: 12.42.tex
\subsection{Symmetricity}
\label{subsection;1.14.20}

We introduce another $\gbigs_n$-action $\rho$
on $H^{\ast}(C^n)[\omega_1,\ldots,\omega_n]$
defined as follows:
\begin{quote}
For $\sigma\in\gbigs_n$,
we put
$\rho(\sigma)\bigl(\prod_i p_i^{\ast}(a_i)\cdot\omega_i^{l_i}\bigr)
  :=\prod_i p_{\sigma(i)}^{\ast}(a_i)\cdot\omega_{\sigma(i)}^{l_i}
 $.
\end{quote}
The action $\rho$ preserves the products.
Let $H^{\ast}(C^n)[\omega_1,\ldots,\omega_n]^{\rho,\gbigs_n}$
denote the $\gbigs_n$-invariant part of with respect to the action
$\rho$.

\begin{thm}\label{thm;1.14.21}
We have $\Psi^{\ast}H^{\ast}(\Quot(\infty,n))=
 H^{\ast}(C^n)[\omega_1,\ldots,\omega_n]^{\rho,\gbigs_n}$.
\end{thm}
\pf
Let consider the element $\xi(l\cdot \vece_n)$.
We have the following relation:
\begin{equation} \label{eq;1.10.1}
 \xi(l\cdot \vece_n)=
 \omega_n\cdot \xi((l-1)\cdot \vece_n)
+\sum_{k<n}\Delta_{n,k}\cdot \xi((l-1)\cdot \vece_k).
\end{equation}
Let $t$ be a formal variable.
We put as follows:
\[
 P_n(t):=\sum_{l=0}^{\infty}t^l\cdot \xi(l\cdot \vece_n).
\]
\begin{lem}
The following holds:
\begin{equation} \label{eq;1.14.26}
 P_n(t)=
 \sum_{h=0}^{n-1}
 \sum_{\substack{J\subset [1,n-1]\\|J|=h}}
 \Delta_{(J,n)}\cdot t^h
 \cdot
 \prod_{i\in J}(1-\omega_i\cdot t)^{-1}
 \cdot
 (1-\omega_n\cdot t)^{-1}.
\end{equation}
Here $(J,n)$ denote the set $J\sqcup \{n\}$.
\end{lem}
\pf
The equality (\ref{eq;1.9.71}) implies the following:
\begin{multline} \label{eq;1.10.2}
P_n(t)=1
+\sum_{l=1}^{\infty}t^l\cdot
 \Bigl(\omega_n\cdot \xi((l-1)\cdot \vece_n)
+\sum_{k<n}\Delta_{n,k}\cdot \xi((l-1)\cdot \vece_k)
 \Bigr)\\
=1+\omega_n\cdot t\cdot P_n(t)
+\sum_{k<n}\Delta_{n,k}\cdot t\cdot P_n(t).
\end{multline}
Namely we have the following:
\begin{equation}\label{eq;1.14.25}
 P_n(t)=
 (1-\omega_n\cdot t)^{-1}+
 \sum_{k<n}\Delta_{n,k}\cdot 
 (1-\omega_n\cdot t)^{-1}\cdot
 t\cdot P_k(t)
\end{equation}
A successive use of the equality (\ref{eq;1.14.25})
provides the formula (\ref{eq;1.14.26}).
\hfill\qed

Then
we have the following relation for any $a\in H^{\ast}(C)$:
\[
 \sum_{n\geq 1}\sum_{l\geq 0}
 t^l\cdot \xi(l\cdot \vece_n)\cdot p_n^{\ast}(a)
=\sum_{h\geq 1}\sum_{|I|=h}
 \Delta_I(a)
 \cdot t^{h-1}\cdot \prod_{j\in I}(1-\omega_j\cdot t)^{-1}.
\]
Here $\Delta_I(a)$ denotes the cohomology class
$\Delta_I\cdot p_i^{\ast}(a)$ for some $i\in I$.
Note that it is independent of a choice of $i$.
Since the right hand side is symmetric,
we obtain that $\sum_{n\geq 0}\xi(l\cdot \vece_n)\cdot p_n^{\ast}(a)$
is symmetric.

Since $S_2=
 \{\xi^Q(l\cdot\vece_1,p^{\ast}a)\,|\,
 a\in H^{\ast}(C),\,\,l=1,\ldots,n\}$
generates the ring $H^{\ast}(\Quot(E,n))$ over $H^{\ast}(C^{(n)})$
(Proposition \ref{prop;1.13.11}),
we obtain the implication
$\Psi^{\ast}H^{\ast}(\Quot(E,n))\subset
 H^{\ast}(C^n)[\omega_1,\ldots,\omega_n]^{\rho,\gbigs_n}$.

Then we only have to prove that
$\Psi^{\ast}S_2$ generates
$H^{\ast}(C^n)[\omega_1,\ldots,\omega_n]^{\rho,\gbigs_n}$.
We only indicate an outline of the argument.

\begin{enumerate}
\item
 For $J=(j_1,\ldots,j_n)\in\seisuuplus^n$,
 we put $\omega^J=\prod_{i=1}^n\omega_i^{j_i}$.
 Any element of
     $H^{\ast}(C^n)[\omega_1,\ldots,\omega_n]^{\rho,\gbigs_n}$
is described as
 $\sum_{J\in\nbigb(n)}\sum_{\sigma\in \gbigs_n}
 \rho(\sigma)(a_J\cdot \omega^J)$.
\item
 For a positive integer $i$ and any element $a\in H^{\ast}(C^{(i)})$,
 we put as follows:
 \[
  e(i,a):=
 \frac{i!(n-i)! }{n!}
 \sum_{\sigma\in\gbigs_n}
 \sigma\Big(a\cdot \prod_{j=1}^i\omega_j\Big).
\]
Then $S_3:=\{ e(i,a)\,|\,i=1,\ldots,n,\,\,a\in H^{\ast}(C^{(i)})\}$
 generates $H^{\ast}(C^n)[\omega_1,\ldots,\omega_n]^{\rho,\gbigs_n}$
 over $H^{\ast}(C^{(n)})$.

It is shown by a similar argument to the proof of Proposition 
 \ref{prop;1.13.10}.

\item
 For a positive integer $i$ and any element $a\in H^{\ast}(C)$,
 we put $p(j,a):=\sum_{i=1}^n p_i^{\ast}(a)\cdot\omega_i^j$.
 Then $S_4:=\{p(j,a)\,|\,j=1,\ldots,n,\,\,a\in H^{\ast}(C)\}$
 generates $H^{\ast}(C^n)[\omega_1,\ldots,\omega_n]^{\rho,\gbigs_n}$
 over $H^{\ast}(C^{(n)})$.

It can be shown by an argument similar to the proof of 
Proposition \ref{prop;1.13.11}.

\item
Let $F^{(3)}_hH^{\ast}(C^{n})[\omega_1,\ldots,\omega_n]^{\rho,\gbigs_n}$
denote the subspace generated
by $\big\{\sum_{J\in \nbigb(n)}\rho(\sigma)(a_j\cdot\omega^J)\,|\,
 \co(J)\leq h\big\}$.
The filtration $F^{(3)}$ is compatible with the product.
The elements $\xi(l\cdot\vece_1,a)$ and $p(l,a)$
are same in the associated graded vector space $\Gr^{(3)}_l$.
Then we can conclude that $S_2$ generates 
$H^{\ast}(C^{n})[\omega_1,\ldots,\omega_n]^{\rho,\gbigs_n}$.
\end{enumerate}
Thus we are done.
\hfill\qed

\subsection{A combinatorial description of $\Psi^{\ast}\xi^Q(\vecv;a)$}
\label{subsection;1.15.1}
\subsubsection{A tuple of of subsets}

\begin{df}\mbox{{}}
\begin{enumerate}
\item
Let $\nbigi=(I_1,\ldots,I_l)$ be a tuple of subsets of $[1,n]$,
namely $I_j\subset[1,n]$.
We say that $I_j$ and $I_k$ are connected in $\nbigi$,
if we have a numbers $i_1,i_2,\ldots,i_h$ such that
$I_j\cap I_{i_1}\neq\emptyset$,
$I_{i_h}\cap I_k\neq \emptyset$,
and
$I_{i_{\alpha}}\cap I_{i_{\alpha+1}}\neq \emptyset$
$(\alpha=1,\ldots,h-1)$.
\item
It determines the equivalence relation on $\nbigi$.
The equivalence class is called the connected component.
The set of the connected components is denoted by $c(\nbigi)$.
Thus we have the decomposition
$\nbigi=\coprod_{\nbigi_{\alpha}\in c(\nbigi)}\nbigi_{\alpha}$.
\hfill\qed
\end{enumerate}
\end{df}

\noindent
{\bf Example}\\
Let consider the case $I_1=\{1,2\}$, $I_2=\{2,3\}$, $I_3=\{4\}$
and $\nbigi=(I_1,I_2,I_3)$.
Then $I_1$ and $I_2$ are connected, and $I_1$ and $I_3$ are not connected.
Hence the connected components are
$\nbigi_1=(I_1,I_2)$ and $\nbigi_2=(I_3)$.

\begin{df}
Assume that $\nbigi=(I_1,\ldots,I_l)$ is connected.
Then we obtain the connected graph $G(\nbigi)$ as follows:
\begin{quote}
The set of vertices are $\{I_1,\ldots,I_l\}\cup \bigcup_{j=1}^l I_j$.
The set of the edges are given as follows:
The vertex $I_i$ and $a\in \bigcup_{j=1}^l I_j$ is connected by an edge
if and only if $a\in I_i$.
\end{quote}
The first betti number of the graph $G(\nbigi)$
is denoted by $b_1(\nbigi)$.
\hfill\qed
\end{df}

\noindent
{\bf Example}
($I_1=\{1,2\}$, $I_2=\{2,3\}$ and $\nbigi=(I_1,I_2)$)
\begin{center}
 \begin{picture}(200,70)(0,-15)
 \put(30,0){\line(-1,1){30}}
 \put(30,0){\line(1,1){30}}
 \put(90,0){\line(-1,1){30}}
 \put(90,0){\line(1,1){30}}
 \put(0,35){\mbox{1}}
 \put(60,35){\mbox{2}}
 \put(120,35){\mbox{3}}
 \put(30,-10){\mbox{$I_1$}}
 \put(90,-10){\mbox{$I_2$}}
 \put(150,20){\mbox{$b_1(\nbigi)=0$}}
 \put(28,-2){\mbox{$\bullet$}}
 \put(88,-2){\mbox{$\bullet$}}
 \put(-2,28){\mbox{$\bullet$}}
 \put(58,28){\mbox{$\bullet$}}
 \put(118,28){\mbox{$\bullet$}}
 \end{picture}
\end{center}

\noindent
{\bf Example}
($I_1=\{1,2\}$, $I_2=\{2,3\}$, $I_3=\{1,3\}$,
 $\nbigi=(I_1,I_2,I_3)$)
\begin{center}
\begin{picture}(200,100)(0,-15)
 \put(30,0){\line(-1,1){30}}
 \put(30,0){\line(1,1){30}}
 \put(90,0){\line(-1,1){30}}
 \put(90,0){\line(1,1){30}}
 \put(60,60){\line(-2,-1){60}}
 \put(60,60){\line(2,-1){60}}
 \put(-5,35){\mbox{1}}
 \put(60,35){\mbox{2}}
 \put(120,35){\mbox{3}}
 \put(30,-10){\mbox{$I_1$}}
 \put(90,-10){\mbox{$I_2$}}
 \put(55,65){\mbox{$I_3$}}
 \put(28,-2){\mbox{$\bullet$}}
 \put(88,-2){\mbox{$\bullet$}}
 \put(-2,28){\mbox{$\bullet$}}
 \put(58,28){\mbox{$\bullet$}}
 \put(118,28){\mbox{$\bullet$}}
 \put(58,58){\mbox{$\bullet$}}
 \put(150,20){\mbox{$b_1(\nbigi)=1$}}
\end{picture}
\end{center}

\noindent
{\bf Example}
($I_1=\{1,2\}$, $I_2=\{2,3\}$, $I_3=\{1,2,3\}$,
 $\nbigi=(I_1,I_2,I_3)$)
\begin{center}
\begin{picture}(200,100)(0,-15)
 \put(30,0){\line(-1,1){30}}
 \put(30,0){\line(1,1){30}}
 \put(90,0){\line(-1,1){30}}
 \put(90,0){\line(1,1){30}}
 \put(60,60){\line(-2,-1){60}}
 \put(60,60){\line(0,-1){30}}
 \put(60,60){\line(2,-1){60}}
 \put(-5,35){\mbox{1}}
 \put(65,35){\mbox{2}}
 \put(120,35){\mbox{3}}
 \put(30,-10){\mbox{$I_1$}}
 \put(90,-10){\mbox{$I_2$}}
 \put(55,65){\mbox{$I_3$}}
  \put(28,-2){\mbox{$\bullet$}}
 \put(88,-2){\mbox{$\bullet$}}
 \put(-2,28){\mbox{$\bullet$}}
 \put(58,28){\mbox{$\bullet$}}
 \put(118,28){\mbox{$\bullet$}}
 \put(58,58){\mbox{$\bullet$}}
 \put(150,20){\mbox{$b_1(\nbigi)=2$}}
\end{picture}
\end{center}

Let $\nbigi=(I_1,\ldots,I_l)$ be a tuple, which is not necessarily
connected.
For each connected component $\nbigi_{\alpha}\in c(\nbigi)$,
we obtain the first betti number $b_1(\nbigi_{\alpha})$.
We put
$c_i(\nbigi):=
 \{\nbigi_{\alpha}\in c(\nbigi)\,|\,
 b_1(\nbigi_{\alpha})=i
 \}$.
Then we have the decomposition
$c(\nbigi)=\coprod_i c_i(\nbigi)$.

\begin{df}
In general,
we put $S(\nbigi):=\bigcup_{I_j\in\nbigi}I_j$,
which is called the support of $\nbigi$.
\hfill\qed
\end{df}

\subsubsection{Combinatorial data}

For an element $\vecl\in\seisuu^n$,
the $j$-th component is denoted by $q_j(\vecl)$.
For an element $\vecl\in\seisuuplus^n$,
we put $s(\vecl):=\{i\,|\,l_i\neq 0\}$.
It is called the support of $\vecl$.
Let consider a tuple
$\vecL=(\vecl_1,\ldots,\vecl_n)$,
where $\vecl_i$ is an element of $\seisuuplus^i$.
For such a tuple, we put $\shat(\vecl_i)=s(\vecl)\cup\{i\}$.
Then we obtain a tuple of subsets of $[1,n]$:
\[
 \nbigi(\vecL):=
 \bigl(
 \shat(\vecl_1),\shat(\vecl_2),\ldots,\shat(\vecl_n)
 \bigr).
\]

We denote the projection of $\seisuuplus^n$
onto the $j$-th component by $q_j$.
\begin{df}
For an element $\vecu\in \nbigb(n)$ and $\sigma\in\gbigs_n$,
we consider the following conditions:
\begin{itemize}
\item
 $\sigma(\sum_{j=1}^n \vecl_j)=\vecu$.
\item
 $c(\nbigi(\vecL))=c_1(\nbigi(\vecL))\sqcup c_2(\nbigi(\vecL))$.
 In other words,
 $c_i(\nbigi(\vecL))=\emptyset$ for any $i\geq 2$.
\item
 We put $\vecl(h)=\sum_{j\leq h}\vecl_j$.
 Let $i$ and $j$ be elements of $\shat(\vecl_h)$,
 and then $q_i(\vecl(h))\neq q_j(\vecl(h))$.
 If $q_j(\vecl(h))<q_i(\vecl(h))$,
 then the inequality $q_j(\vecl(h))\leq q_i(\vecl(h-1))$ holds.
\end{itemize}
The set of
$\vecL=(\vecl_1,\ldots,\vecl_n)$
satisfying the conditions
is denoted by $\nbigt(\vecu,\sigma)$.
\hfill\qed
\end{df}

For an element $\vecL\in \nbigt(\vecu,\sigma)$,
we put $\rho_h(\vecL):=|\vecl_h|-|\shat(\vecl_h)|+1$.

\subsubsection{A formula}

\begin{thm}\label{thm;1.14.45}
The following holds:
\begin{equation}\label{eq;1.14.46}
 \Psi^{\ast}\xi^Q(\vecu;a)
=\frac{1}{|St(\vecu)|}
 \sum_{\sigma\in\gbigs_n}
 \sum_{\vecL\in\nbigt(\vecu,\sigma)}
 \prod_{h=1}^n\omega_h^{\rho_h(\vecL)}\cdot
 \prod_{\nbigi_0\in c_0(\vecL)}\Delta_{\nbigi_0}
\cdot
 \prod_{\nbigi_1\in c_1(\vecL)}
 \big((2-2g)\cdot pt_{\nbigi_1}\bigr)\cdot\sigma(a).
\end{equation}
Here $g$ denotes the genus of $C$.
\end{thm}
\pf
The proof is just a calculation.
We only indicate an outline in the next two subsubsections.
\subsubsection{A reduction}

\label{subsubsection;1.15.20}

Let $\vecu$ be an element of $\seisuuplus^{n-1}$.
We have the following formula:
\begin{multline}\label{eq;1.14.30}
 \xi(\vecu+l\cdot\vece_n)\cdot a=
 \omega_n\cdot\xi(\vecu+(l-1)\cdot\vece_n)\cdot a
+\sum_{\substack{k<n\\u_k\leq l-1}}
 \Delta_{n,k}\cdot\xi\bigl(\tau_{n,k}(\vecu+(l-1)\cdot\vece_n)\bigr)\cdot a\\
=\omega_n\cdot\xi(\vecu+(l-1)\cdot\vece_n)\cdot a
+\sum_{\substack{k<n\\u_k\leq l-1}}
 \Delta_{n,k}\cdot\xi\bigl(\tau_{n,k}(\vecu+(l-1)\cdot\vece_n)\bigr)\cdot
 \tau_{n,k}(a).
\end{multline}
Here we use $\Delta_{n,k}\cdot a=\Delta_{n,k}\tau_{n,k}(a)$.

Let $\nbigp$ be a tuple $(k_1,\ldots,k_h;l_1,\ldots,l_h)$
satisfying the following:
\[
 u_{k_1}+l_1<l,\,\,
 u_{k_2}+l_2<u_{k_1},\,
\ldots,
 u_{k_h}+l_h<u_{k_h-1}.
\]
Let $\nbigu$ denote the set of such tuple.
For such a tuple, we put $|\nbigp|:=\sum_{j=1}^hl_j+u_{k_h}-h$,
and $\Delta_{\nbigp}=\Delta_{(k_1,\ldots,k_h,n)}$.
The cyclic permutation $\sigma_{\nbigp}$
is defined by 
$n\mapsto k_1$,
$k_i\mapsto k_{i+1}$
and $k_n\mapsto n$.
We put as follows:
\[
 \psi_{\nbigp}(\vecu+l\cdot\vece_n)=
 \sigma_{\nbigp}
 \Big(\vecu+(l-l_n)\cdot\vece_n-
 \sum_{j=1}^{h-1}l_{j+1}\vece_{k_j}
 -u_{k_n}\cdot\vece_{k_n}
 \Big).
\]
\begin{lem}
The following holds:
\[
 \xi(\vecu+l\cdot \vece_n)\cdot a
=\sum_{\nbigp\in\nbigu}
 \omega_n^{|\nbigp|} \cdot
 \Delta_{\nbigp}\cdot
 \xi(\psi_{\nbigp}(\vecu+l\cdot\vece_n))\cdot\sigma_{\nbigp}(a).
\]
\end{lem}
\pf
We only have to use the formula (\ref{eq;1.14.30})
recursively.
\hfill\qed

\subsubsection{Statement for an induction}

Let $h$ be an integer such that $0\leq h\leq n$.
Let $I=(i_{n-h+1},i_{n-h+2},\ldots,i_n)$ be a tuple such that
$i_{\alpha}\in [1,n]$.
We put $\gbigs(n,I):=
 \{\sigma\in\gbigs_n\,|\,\sigma(i_{\alpha})=\alpha\}$. 
Then we have the decomposition
$\gbigs_n=\coprod_{|I|=h}\gbigs(n,I)$.

\begin{df}
Let $h$ be an integer such that $0\leq h\leq n$.
Let $\vecv$ be an element of $\nbigb(n)$
and $\sigma$ be an element of $\gbigs_n$.

Let consider a tuple
$\vecL=(\vecl_{n-h+1},\vecl_{n-h+2},\ldots,\vecl_n)$
such that $\vecl_{\alpha}\in \seisuuplus^{\alpha}$.
Let consider the following conditions:
\begin{itemize}
\item
 $\sum_{\alpha=j}^nq_j(\vech_{\alpha})=
 q_j(\sigma(\vecv))$
for any $j\geq n-h+1$.
In other words,
 $\sigma(\vecv)-\sum_{\alpha=n-h+1}^n\vech_{\alpha}\in\seisuuplus^{n-h}$.
\item
 We put $\shat(\vech_{\alpha})=s(\vech_{\alpha})\cup \{\alpha\}$.
 We put $\vecv^{(\beta)}
=\sigma(\vecv)-\sum_{\alpha=\beta}^n\vech_{\alpha}$.
 Let $i$ and $j$ be elements of $\shat(\vech_{\alpha})$.
 Then $q_j(\vecv^{(\alpha)})\neq q_i(\vecv^{(\alpha)})$.
 If $q_j(\vecv^{(\alpha)})<q_i(\vecv^{(\alpha)})$,
then the inequality $q_j(\vecv^{(\alpha)})\leq q_i(\vecv^{(\alpha-1)})$
holds.
\end{itemize}
The set of such tuples $\vecL$ is denoted by
$\nbigt_h(\vecv,\sigma)$.
\hfill\qed
\end{df}
For such a tuple $\vecL=(\vecl_{n-h+1},\ldots,\vecl_n)$,
we put
$\rho_{\alpha}(\vecL)
:=\sum_{i=1}^nq_i(\vech_{\alpha})
 -|\shat(\vech_{\alpha})|+1$,
$(\alpha=n-h+1,\ldots,n)$.

For any element $\vecv\in\nbigb(n)$ and $a\in H^{\ast}(F^Q(\vecv))$,
we put as follows:
\[
 \xi^{\sym}(\vecv;a):=
 \sum_{\sigma\in\gbigs_n}
 \xi(\sigma(\vecv))\cdot \sigma(a).
\]

\begin{lem}
Let $h$ be an integer such that $0\leq h\leq n$.
Then we have the following:
\begin{equation}
 \xi^{\sym}(\vecv;a)=
\sum_{|I|=h}
\sum_{\sigma\in \gbigs(n,I)}
\sum_{\vecL\in\nbigt_h(\vecv,\sigma)}
 \prod_{\alpha=n-h+1}^n
 \Bigl(
 \omega_{\alpha}^{\rho_{\alpha}(\vecL)}
 \cdot
 \Delta_{\shat(\vech_{\alpha})}
 \Bigr)
 \cdot\xi\Big(\sigma(\vecv)-\sum_{\alpha=n-h+1}^n\vech_{\alpha}\Big)
 \cdot\sigma(a).
\end{equation}
\end{lem}
\pf
We use an ascending induction on $h$.
We only have to use a reduction given
in the subsubsection \ref{subsubsection;1.15.20}
recursively.
\hfill\qed

In the case $h=n$,
the set $\gbigs(\sigma,I)$ contains the unique element,
and we have $\nbigt_n(\vecu,\sigma)=\nbigt(\vecu,\sigma)$.
Thus we obtain the following equality:
\begin{equation}\label{eq;1.14.40}
 \xi^{\sym}(\vecv;a)=
\sum_{\sigma\in\gbigs_n}
\sum_{\vecL\in\nbigt(\vecv,\sigma)}
 \prod_{\alpha=1}^n
 \Bigl(
 \omega_{\alpha}^{\rho_{\alpha}(\vecL)}
 \cdot
 \Delta_{\shat(\vech_{\alpha})}
 \Bigr)
 \cdot\xi\Big(\sigma(\vecv)-\sum_{\alpha=1}^n\vech_{\alpha}\Big)
 \cdot\sigma(a).
\end{equation}
The following lemma is easy.
\begin{lem}
Let $\nbigi=(I_1,\ldots,I_h)$ be a connected tuple
of subsets of $[1,n]$.
The following holds:
\begin{equation}\label{eq;1.14.41}
 \prod_{\alpha=1}^a
 \Delta_{I_{\alpha}}
=
\left\{
\begin{array}{ll}
 0, & (b_1(\nbigi)\geq 2),\\
 (2-2g)\cdot pt_{S(\nbigi)}, & (b_1(\nbigi)=1),\\
 \Delta_{S(\nbigi)}, & (b_1(\nbigi)=0).
\end{array}
\right.
\end{equation}
Here $g$ denotes the genus of $C$.
\hfill\qed
\end{lem}

By substiting (\ref{eq;1.14.41}) into
(\ref{eq;1.14.40}),
we obtain (\ref{eq;1.14.46}).
Thus Theorem \ref{thm;1.14.45} is obtained.
\hfill\qed

%% file: 12.41.tex
\subsection{General case}

For an element $\vecl=(l_1,\ldots,l_h)\in\seisuuplus^h$,
we put $\gbigs(\vecl):=\prod_{i=1}^h\gbigs_{l_i}$.
We have the natural action of $\gbigs(\vecl)$
on $\prod_{i=1}^h\seisuuplus^{l_i}$.
For an element
$\vecv_{\ast}=(\vecv_1,\ldots\vecv_h)\in
\prod_{i=1}^h\seisuuplus^{l_i}$,
we put $S(\vecv_{\ast})=\prod_{i=1}^h St(\vecv_i)$.

We naturally regard $\gbigs(\vecl)$
as the subgroup of $\gbigs(|\vecl|)$.
We also naturally identify
$\prod_{i=1}^{h}\seisuuplus^{l_i}$
with $\seisuuplus^{|\vecl|}$.

For an element $\vecl=(l_1,\ldots,l_h)$,
we put $\rho_j(\vecl)=\sum_{i\leq j}l_i$.
We have $\rho_h(\vecl)=|\vecl|$.
We have the natural morphism
$\Psi_{\vecl}:\Filt(E,\vecl)\lrarr \Filt(E,|\vecl|)$,
given by the following functor:
\[
 \Bigl(
 E\rarr
 \nbigf_{|\vecl|}\rarr\nbigf_{|\vecl|-1}\rarr\cdots\rarr \nbigf_1
 \Bigr)
\longmapsto
 \Bigl(
 E\rarr \nbigf_{\rho_h(\vecl)}\rarr\nbigf_{\rho_{h-1}(\vecl)}
 \rarr\cdots\rarr\nbigf_{\rho_1(\vecl)}
 \Bigr).
\]
Thus we obtain the morphism
$\Psi_{\vecl}^{\ast}:
 H^{\ast}(\Filt(E,\vecl))\lrarr
 H^{\ast}(\Filt(E,|\vecl|))$.

\begin{prop}\mbox{{}}
\begin{itemize}
\item
For any element $\vecv_{\ast}\in\nbigb(\vecl,r)$
and for any $a\in H^{\ast}(F^{\vecl}(\vecv_{\ast}))$,
we have the following equality:
\[
 \Psi_{\vecl}^{\ast}
 \xi^{\vecl}(\vecv_{\ast};a)
=\frac{1}{|St(\vecv_{\ast})|}
\sum_{\sigma\in\gbigs(\vecl)}
 \sigma \bigl(
 \xi(\vecv_{\ast})\cdot a\bigr).
\]
\item
The filtration
$\left\{F_hH^{\ast}(\Filt(E,\vecl))\,|\,h=0,1,2,\ldots,\right\}$
is compatible with the product,
that is,
$F_h\cdot F_k\subset F_{h+k}$.

\item
Thus we obtain the associated graded ring
$Gr H^{\ast} (\Filt(E,\vecl))$.
It is naturally identified with
the $\gbigs(\vecl)$-invariant part
of $Gr H^{\ast}(\Filt(E,|\vecl|))$.
\hfill\qed
\end{itemize}
\end{prop}

On $\Filt(E,\vecl)\times C$,
we have the universal quotients
$q_2^{\ast}E\rarr
\nbigfu_{h}\rarr\nbigfu_{h-1}\rarr\cdots\rarr\nbigfu_{1}$.
We put $\nbiggu_j:=\ker(\nbigfu_j\lrarr\nbigfu_{j-1})$.
We have the $i$-th Chern class
$c_i(\nbiggu_j)$.
It induces the correspondence map
$\phi_i^{(j)}:H^{\ast}(C)[-2(i-1)]\lrarr H^{\ast}(\Filt(E,\vecl))$.
\begin{thm}
The set 
$ \bigcup_{j=1}^n\bigcup_{i=1}^{l_j+1}
 \Image \phi_i^{(j)}$
generates $H^{\ast}(\Filt(E,\vecl))$.
\hfill\qed
\end{thm}

The $\gbigs_{|\vecl|}$-action $\rho$ on
$H^{\ast}(C^{|\vecl|})[\omega_1,\ldots,\omega_{|\vecl|}]$
induces the $\gbigs(\vecl)$-action,
which we denote by $\rho$.
We denote
the $\gbigs(\vecl)$-invariant part with respect to $\rho$
by 
$H^{\ast}(C^{|\vecl|})[\omega_1,\ldots,\omega_{|\vecl|}]
 ^{\rho,\gbigs(\vecl)}$.
\begin{thm}
The image $\Psi_{\vecl}^{\ast}H^{\ast}(\Filt(E,\vecl))$
is $H^{\ast}(C^{|\vecl|})[\omega_1,\ldots,\omega_{|\vecl|}]
 ^{\rho,\gbigs(\vecl)}$.
\hfill\qed
\end{thm}

We omit to derive a combinatorial description
of $\Psi_{\vecl}^{\ast}\xi^{\vecl}(\vecv_{\ast},a)$.

%% file: 12.6.tex
\section{The structure of the cohomology ring of
the infinite quot scheme of infinite length}

Recall the following general notation:
Let $V^{\cdot}$ be a graded vector space such that
$V^i=0$ for any $i<0$ and that
$\dim V^i<\infty$ for any $i$.
Recall that the formal power series $\sum_{i\geq 0} (\dim V^i)\cdot t^i$
is called the Poincar\'e series for $V^{\cdot}$.
Let $X$ be a topological space whose $i$-th betti number is finite
for any $i$.
We denote the Poincar\'e series of the 
graded vector space $H^{\ast}(X)$ by $P(X)(t)$,
or simply by $P(X)$.

Let $C$ be a smooth projective curve.
Fix a point $P$ of $C$.
We have the inclusion of $\Quot(r,r)$ into $\Quot(r+1,r+1)$
corresponding to the functor:
\[
 \left(\nbigo^{\oplus r}\lrarr \nbigf\right)
\longmapsto
\left(
 \left(\nbigo^{\oplus r}\lrarr \nbigf\right)
\oplus
 (\nbigo\lrarr \cnum_P)
\right),
\]
Here $\cnum_P$ denotes the skyscraper sheaf at the point $P$.
Thus we obtain the topological space
$\Quot(\infty,\infty):=\bigcup_{r=1}^{\infty}\Quot(r,r)$,
which we call the infinite quot scheme of infinite length
for $C$.
(More precisely, we use the mapping telescope.)
We consider the structure of the cohomology ring
of $\Quot(\infty,\infty)$.

We put $b_1:=\dim H^1(C)$.
\begin{lem}
The Poincar\'e series of $\Quot(\infty,\infty)$ is the following,
\[
 \frac{(1+t)^{b_1}}{(1-t^2)}
 \prod_{h=1}^{\infty}
\frac{(1+t^{2h+1})^{b_1}}{(1-t^{2h})(1-t^{2h+2})}
\in \seisuu[t][[s]].
\]
\end{lem}
\pf
We have the isomorphism
$H^{\ast}(\Quot(r,l))\simeq 
 \bigoplus_{\vecl\in \Dec(l,r)}
  H^{\ast}(C^{(\vecl)})[-co(\vecl)]$.
Then we obtain the following equality:
\[
 \sum_{l=0}^{\infty}P\big(Quot(r,l)\big)\cdot s^l=
 \prod_{h=0}^{r-1}
 \frac{(1+st^{2h+1})^{b_1}}{(1-t^{2h}s)(1-t^{2h+2}s)}
\in\seisuu[t][[s]].
\]
And thus we have the following equality:
\[
 \sum_{l=0}^{\infty}P\big(Quot(\infty,l)\big)\cdot s^l=
 \prod_{h=0}^{\infty}
 \frac{(1+st^{2h+1})^{b_1}}{(1-t^{2h}s)(1-t^{2h+2}s)}
\in \seisuu[[t]][[s]].
\]
The following holds:
\[
 P(\Quot(\infty,\infty))=
 \sum_{r=0}^{\infty}
 \Bigl(P\big(\Quot(r+1,r+1)\big)-P\big(\Quot(r,r)\big)\Bigr).
\]
Here  $\Quot(0,0)=\emptyset$.
Since we have
${\displaystyle
 \lim_{r\to\infty}
  \Bigl(P\big(\Quot(\infty,r)\big)-P\big(\Quot(r,r)\big)\Bigr)=0
 }$
in the ring $\seisuu[[t]]$,
the right hand side is same as the following:
\[
 \sum_{r=0}^{\infty}
 \Bigl(P\big(\Quot(\infty,r+1)\big)-P\big(\Quot(\infty,r)\big)
 \Bigr).
\]
Here $\Quot(\infty,0)=\emptyset$.
Note that the following holds:
\[
 (1-s)\sum_{l=0}^{\infty}P\big(Quot(\infty,l)\big)\cdot s^l=
\sum_{l=1}^{\infty}
\Bigl( P\big(\Quot(\infty,l)\big)
  -P\big(\Quot(\infty,l-1)\big) \Bigr)\cdot s^l.
\]
We have the following equalities,
\begin{multline}
 (1-s)\sum _{l=0}^{\infty}P\big(\Quot(\infty,l)\big)\cdot s^l=
 (1-s)\cdot\frac{(1+st)^{b_1}}{(1-s)(1-st^2)}
 \prod_{h=1}^{\infty}
 \frac{(1+st^{2h+1})^{b_1}}{(1-t^{2h}s)(1-t^{2h+2}s)}\\
=\frac{(1+st)^{b_1}}{(1-st^2)}
 \prod_{h=1}^{\infty}
 \frac{(1+st^{2h+1})^{b_1}}{(1-t^{2h}s)(1-t^{2h+2}s)}. 
\end{multline}
It is contained in $\seisuu[s][[t]]$.
Hence we can specialize $s=1$.
Thus we obtain the following equality:
\begin{multline}
P\big(\Quot(\infty,\infty)\big)=
\Bigl.(1-s)\sum_{l=0}^{\infty}
   P\big(\Quot(\infty,l)\big)\cdot s^l\Bigr|_{s=1}\\=
\left.\frac{(1+st)^{b_1}}{(1-st^2)}
 \prod_{h=1}^{\infty}
 \frac{(1+st^{2h+1})^{b_1}}{(1-t^{2h}s)(1-t^{2h+2}s)}\right|_{s=1}
\!\!
=\frac{(1+t)^{b_1}}{(1-t^2)}
 \prod_{h=1}^{\infty}
\frac{(1+t^{2h+1})^{b_1}}{(1-t^{2h})(1-t^{2h+2})}.
\end{multline}
Thus we are done.
\hfill\qed

\begin{thm}\label{thm;5.18.4}
The rational cohomology ring $H^{\ast}(\Quot(\infty,\infty))$ is
isomorphic to the following,
\[
 H^{\ast}(C^{(\infty)})
\otimes
 \bigotimes_{h=1}^{\infty}
 \Sym^{\cdot}\bigl(H^{\ast}(C)[-2h]\bigr).
\]
\end{thm}
\pf
The natural morphism $\Quot(r,r)\lrarr C^{(r)}$ induces
the morphism:
\[
 H^{\ast}(C^{(\infty)})\lrarr H^{\ast}(\Quot(\infty,\infty)).
\]
We have the morphism
$H^{\ast}(C)[-2h]\lrarr H^{\ast}(\Quot(r,r)),
a\longmapsto \xi(h\cdot \vece_1;a)$.
Thus we obtain the ring morphism:
\[ 
 H^{\ast}(C^{(\infty)})\otimes \bigotimes_{h=1}^{\infty}
 Sym^{\cdot}\bigl(H^{\ast}(C)[-2]\bigr)
\lrarr 
  H^{\ast}(\Quot(\infty,\infty)).
\]
Thus we obtain the ring morphism:
\[ 
 H^{\ast}(C^{(\infty)})\otimes
 \bigotimes_{h=1}^{\infty}
 \Sym^{\cdot}\bigl(H^{\ast}(C)[-2]\bigr)
\lrarr 
  H^{\ast}(Quot(\infty,\infty)).
\]
Due to the result in the previous section,
it is surjective.
The Poincar\'e series of the both sides coincide.
Thus we are done.
\hfill\qed

%% file: quot3_new.bbl
\begin{thebibliography}{99}

\bibitem{bf}
K. Behrend and B. Fantechi,
{\em The intrinsic normal cone},
Invent. Math. {\bf 128}, pp 45--88,
(1997).
\bibitem{bb}
A. Bianicki-Birula,
{\it Some theorems on actions of algebraic groups},
Annals of Math. {\bf 98} (1973),
pp 480--497.
\bibitem{cs}
J. B. Carrel and A. J. Sommes,
{\em Some topological aspects of $\cnum^{\ast}$-action
on compact Kahler manifolds},
Comment, Math. Helv. {\bf 54} (1979),
pp 567--582.


\bibitem{ful} W. Fulton, {\em Intersection theory}
             (Springer, 1984)

\bibitem{gr}
A. Grothendieck,
{\em Techniques de construction et th\'eor\`ems
d'existence en g\'eom\'etrie alg\'ebrique IV:
Les sch\'emas de Hilbert}, S\'eminaire Bourbaki expos\'e 221
(1961), IHP, Paris.
\bibitem{go}
L. G\"ottsche,
{\em Hilbert schemes of zero-dimensional subschemes
of smooth varieties},
Lecture Notes in Mathematics 1572, Springer (1994).


\bibitem{m1}
 T. Mochizuki,
{\em
The cohomology ring of the Quot scheme of coherent sheaves with length 2
over smooth projective curves},
http://www.math.ias.edu/$\,\tilde{}\,$takuro/list.html

\bibitem{m2}
T. Mochizuki,
{\em The cohomology ring of the Quot scheme of coherent sheaves
with length 3 over smooth projective curves},
http://www.math.ias.edu/$\,\tilde{}\,$takuro/list.html


\bibitem{m4}
T. Mochiuzki,
{\em The theory of the invariants obtained from the moduli stacks
of stable objects on a smooth polarized surface},
http://www.math.ias.edu/$\,\tilde{}\,$takuro/list.html

\bibitem{n}
H. Nakajima,
{\em Lectures on Hilbert schemes of points on surfaces},
University Lecture Series, 18.
American Mathematical Society Providence RI, 1999.
\end{thebibliography}
